\documentclass[12pt]{amsart}
\usepackage{amsthm}
\usepackage{amssymb}
\usepackage{MnSymbol}
\usepackage{amsxtra}     
\usepackage{epsfig}
\usepackage{verbatim}
\usepackage[all]{xypic}
\usepackage{url}

\theoremstyle{plain}
\newtheorem{thm}{Theorem}
\newtheorem{prop}[thm]{Proposition}
\newtheorem{cor}[thm]{Corollary}
\newtheorem{lem}[thm]{Lemma}

\newtheorem{warning}[thm]{Warning}
\theoremstyle{definition}
\newtheorem{condition}[thm]{Condition}
\newtheorem{example}[thm]{Example}

\theoremstyle{remark}
\newtheorem{rmk}[thm]{Remark}

\include{header}

\newcommand{\adeles}{\mathbb{A}}

\newcommand{\idelesK}{\mathbb{A}_K^{\times}}
\newcommand{\Ind}{\mbox{\bf Ind}}

\newcommand{\hn}{\mathcal{H}_n}

\newcommand{\OK}{\mathcal{O}_{\cmfield}}
\newcommand{\Oe}{\mathcal{O}_E}
\newcommand{\Hom}{\mbox{Hom}}
\newcommand{\dual}{^\vee}
\newcommand{\isomto}{\overset{\sim}{\rightarrow}}
\newcommand{\lp}{\left(}
\newcommand{\rp}{\right)}
\newcommand{\ci}{C^{\infty}}

\newcommand{\Split}{\mathrm{Split}}

\newcommand{\uo}{\underline{\omega}}

\newcommand{\padic}{p\mathrm{-adic}}

\newcommand{\cmfield}{K}

\newcommand{\uA}{\underline{A}}

\newcommand{\gln}{\mathrm{GL_n}}
\newcommand{\gl}{\mathrm{GL}}

\newcommand{\e}{{\mathbf e}}

\newcommand{\IR}{\mathbb{R}}
\newcommand{\ZZ}{\mathbb{Z}}
\newcommand{\IC}{\mathbb{C}}
\newcommand{\IQ}{\mathbb{Q}}

\newcommand{\End}{\mbox{End}}

\newcommand{\Cp}{\IC_p}
\newcommand{\OCp}{\mathcal{O}_{\Cp}}

\newcommand{\Gal}{\mathrm{Gal}}
\newcommand{\diag}{\mathrm{diag}}
\newcommand{\Aut}{\mathrm{Aut}}
\newcommand{\ord}{\mathrm{ord}}
\newcommand{\cond}{\mathrm{conductor}}
\newcommand{\hern}{\mathrm{Her}_n}
\newcommand{\Sgn}{\mathrm{Sign}}

\renewcommand{\paragraph}[1]{\noindent\textbf{#1}. }
\renewcommand{\subparagraph}[1]{\noindent\emph{#1}.}

\def\tr{{\rm tr}\,}

\begin{document}

\bibliographystyle{amsalpha}       
\title{A $p$-adic Eisenstein Measure for Unitary Groups}
\author{Ellen Eischen}
\maketitle\let\thefootnote\relax\footnote{Revisions to this paper were made while the author was partially supported by National Science Foundation Grant DMS-1249384.}

\begin{abstract}
We construct a $p$-adic Eisenstein measure with values in the space of $p$-adic automorphic forms on certain unitary groups.  Using this measure, we $p$-adically interpolate certain special values of both holomorphic and non-holomorphic Eisenstein series, as both the archimedean and the $p$-adic weights of the Eisenstein series vary.  
\end{abstract}

\setcounter{tocdepth}{2}
\tableofcontents
\setcounter{secnumdepth}{3}

\section{Introduction}
\subsection{Motivating goals}

\subsubsection{Number theory}
This paper completes one step toward a construction of certain $p$-adic $L$-functions on unitary groups.  

Starting with the work of Jean Pierre Serre and Nicholas Katz, $p$-adic familes of Eisenstein series have been used to construct certain $p$-adic $L$-functions.  In \cite{kaCM}, Katz constructs a $p$-adic Eisenstein measure for Hilbert modular forms, which he uses to construct $p$-adic $L$-functions for CM fields.  

Continuing in this direction, we construct a $p$-adic Eisenstein measure with values in the space of $p$-adic automorphic forms on particular unitary groups.  Like in \cite{kaCM}, we $p$-adically interpolate certain special values of both holomorphic and non-holomorphic Eisenstein series (similar to the Eisenstein series studied by Shimura in \cite{sh}), and we allow both the $p$-adic and the archimedean weights to vary.  Unlike in \cite{kaCM} but following the theme of \cite{SHL}, we work with Eisenstein series on adelic groups.  The reasons working adelically are:
\begin{enumerate}
\item{An ongoing joint project of Michael Harris, Jian-Shu Li, Christopher Skinner and myself to construct $p$-adic $L$-functions uses the doubling method, which requires integrating Eisenstein series on adelic groups \cite{EEHLS}.  }
\item{$p$-adic interpolation of special values of the Eisenstein series on adelic unitary groups studied extensively by Goro Shimura (e.g. in \cite{sh}) constitutes a necessary step in the approach of \cite{SHL, EEHLS}.}
\end{enumerate}
In \cite{apptoSHLvv}, we generalize the measure in this paper to the case of Eisenstein series with non-scalar weights.  However, the results in the current paper are sufficient for allowing the variable conventionally denoted $``s"$ in the $L$-function to vary $p$-adically.  They also are sufficient for the conjectured connection with homotopy theory, described below.

\subsubsection{Homotopy theory}
Katz's $p$-adic Eisenstein series for modular forms are one of the main tools in Matthew Ando, Michael Hopkins, and Charles Rezk's study of the Witten genus, a modular form-valued invariant of a particular class of manifolds \cite{hopkins94, hopkinsICM, AHR}.  This invariant gives an orientation in the cohomology theory {\it topological modular forms}.  (The Witten genus was initially defined by Edward Witten in terms of a certain formal power series and was later studied in terms of $p$-adic Eisenstein series by Hopkins.)

Mark Behrens, Hopkins, and Niko Naumann aim to generalize the Witten genus to {\it topological automorphic forms} \cite{beh}.  They expect $p$-adic Eisenstein series on unitary groups, and in particular $p$-adic interpolation of the weight at $\infty$, to be a key ingredient in their work.  The measure constructed here gives the first example of an Eisenstein measure on unitary groups in which the weight at $\infty$ varies $p$-adically.  Although the Eisenstein series in this paper are on $U(n,n)$, they can be pulled back in a natural way to $U(n-1, 1)$, which is the group currently of interest in homotopy theory.

\subsection{Main results}

The main result of this paper is the construction in Section \ref{EisensteinMeasureSection} of a $p$-adic measure with values in the space of $p$-adic automorphic forms on the unitary groups $U(n,n)$.  As part of the construction of this measure, we construct a $p$-adic family of Eisenstein series.  As a corollary, we obtain a $p$-adic measure that $p$-adically interpolates certain special values of Eisenstein series (normalized by an appropriate period) initially defined over $\IC$, including non-holomorphic Eisenstein series.

In Sections \ref{HoloESeriesSection} and \ref{diffopsnew}, we introduce the Eisenstein series with which we work.  In the $\ci$-case, these are similar to the Eisenstein series on adelic groups studied by Shimura in \cite{sh}.  The main change to Shimura's approach, like in \cite{SHL}, is at primes dividing $p$.  As explained below, it was necessary to modify the construction at $p$ in \cite{SHL}.  Note that Lemma \ref{pfflemma} of the current paper provides a formulaic approach to constructing a whole class of Siegel sections at $p$ that are well-suited to $p$-adic interpolation of the Fourier coefficients of the Eisenstein series.

Our approach to $p$-adic Eisenstein series and construction of the $p$-adic Eisenstein measure is a natural generalization of Katz's methods in \cite{kaCM}.  Although the formulas in our paper look more complicated than the formulas derived by Katz in \cite{kaCM}, the reader can plug in $n=1$ to verify easily that we really are studying similar Eisenstein series to those in \cite{kaCM} in that case.  (Indeed, a fair amount of cancellation occurs in the terms in the Fourier coefficients precisely when $n=1$, and the indexing is also simpler in that case.)

\subsection{Relationship with prior results}
The current paper extends, clarifies, and corrects portions of the construction of the measure in \cite{SHL}.  Unlike the construction of \cite{SHL}, the measure in the current paper allows
\begin{enumerate}
\item{$p$-adic interpolation of special values of both holomorphic and {non}-holomorphic Eisenstein series, as the weight of the Eisenstein series varies $p$-adically.}
\item{$p$-adically varying the weight of the Eisenstein series at $\infty$ in addition to $p$.}
\end{enumerate}
These were not goals of the construction in \cite{SHL}.  These generalizations rely on the differential operators in \cite{EDiffOps} and a close adaptation of the approach to Hecke characters in \cite[Section 5]{kaCM}.
\subsubsection{Clarifications and corrections to \cite{SHL}}\label{CorrectionsClarifications}
During the course of this work, I came across several issues in the work in \cite{SHL}.  At Skinner's recommendation, I take this opportunity to give corrections and clarifications to the following points.  

\begin{enumerate}
\item{This paper provides a few corrections to the Fourier coefficients in \cite{SHL}, which are computed as a product of local Fourier coefficients in both \cite{SHL} and the current paper.  In the notation of \cite{SHL}, the most important changes are:
\begin{enumerate}
\item{At archimedean primes, the term $\det\beta^{[E:\IQ]}$ in the local Fourier coefficient in \cite[Equation (3.3.5.3)]{SHL} should be $N_{E/\IQ}(\beta)$.  (Shimura's computation of the local Fourier coefficents shows that each of those terms is of the form $\sigma(\det\beta)$, where $\sigma$ varies over all the embeddings of the totally real field $E$ into $\bar{\IQ}$ \cite{sh, shconfluent}.)  We address this change in Warning \ref{betawarning} and Section \ref{sectionsinfty}.
}

\item{At primes $v$ dividing $p$, the Siegel section in \cite[Equation (3.3.4.7)]{SHL} appears not to be a section of the induced representation $\Ind_{P_v}^{G_v}\left(\chi\left|\cdot\right|^{-s}\right)$, based on evaluation of this section at elements of $P_v$.  (It is possible there is a minor modification that would make it into a section of the induced representation.)  

Moreover, $p$-adically interpolating the non-holomorphic Eisenstein series requires working with a $p$-adic section that records information about determinants that the section in \cite{SHL} loses.  More precisely, we need that the Schwartz function $\phi_\nu$ (defined in \cite[Equation (3.3.4.2)]{SHL}) in terms of which the Siegel section at $p$ is defined has the following property: If $\nu = \left(\det^d, \ldots, \det^d\right),$ then $\phi_{\nu} =\det^d$ on the domain of $\phi_\nu$ in the definite case; this condition matters only for the $p$-adic interpolation, not the definition of the Eisenstein series.  This information about determinants is lost if one uses the definition given in \cite{SHL}.

It also turns out that these changes are necessary in order to carry out the computations of the local Euler factors at $p$ of the $p$-adic $L$-functions in \cite{EEHLS}.

Section \ref{sectionatp} of the current paper provides a class of $p$-adic sections that satisfies these further properties.
}\label{padiccorrections}
\item{Away from $p$ and $\infty$, some of the Fourier coefficients are in terms of polynomials that actually vary as the Hecke character varies.  However, \cite{SHL} seems to treat the values of these polynomials as invariant under $p$-adic variation of the Hecke character.  (See Section \ref{sectionsnotpinfty} of the current paper for formulas for the Fourier coefficients at these places, similar to those given in \cite{SHL}.)  This is closely related to the next item, which concerns variation of the Hecke character.}
\end{enumerate}
}
\item\label{item2}{Let $\chi$ be a Hecke character.  The approach in \cite[Section (3.5)]{SHL} seems to be to fix the archimedean component of a Hecke character on a CM field $K$ and vary just the $p$-adic component.  However, changing only the $p$-adic component of a Hecke character might give a character of the id\`ele group $\idelesK$ that is no longer trivial on $K^\times$.  Instead, we use a similar approach to \cite{kaCM}, which is also necessary for varying the archimedean component of the Eisenstein series in this paper.  

Our method, which is based on the approach in \cite[Sections 5.5, 5.0]{kaCM}, seems only to allow for a construction in which both the $p$-adic and the archimedean components of the weights of the Eisenstein series can vary $p$-adically.  See Sections \ref{viewedpadically} and \ref{EisensteinMeasureSection} for further details on our approach.}

\item{As a consequence of \eqref{item2}, the measure in \cite[Section (3.5)]{SHL} should have been constructed on the product of the group denoted $T(l)$ with the {\it quotient} $(\OK\otimes\ZZ_p)/\overline{\OK^\times}$, where $\overline{\OK^\times}$ denotes the $p$-adic closure of $\OK^\times$.}

\item{
As in the convention for modular forms, our $q$-expansions have coefficients that agree over $\IC$ with the Fourier coefficients of a function of a complex variable in a hermitian symmetric space (e.g. the upper half plane and its generalizations).  This requires normalizing the Eisenstein series on adelic groups by an automorphy factor.  This normalization is missing in \cite{SHL}, which accounts for some of the difference between the $q$-expansion coefficients used in \cite[Theorem (3.5.1)]{SHL} and in Equation \eqref{Gknuchimu} here.  (The normalization in this paper is necessary for applying the $q$-expansion principle to the $q$-expansions, which is necessary for $p$-adically interpolating the Eisenstein series.)

}
\end{enumerate}

\subsection{Acknowledgments}
I am grateful to Matthew Emerton and Christopher Skinner for insightful suggestions during the past year, especially as I came across modifications that I needed to make to \cite{SHL}.  Many of the suggestions found their way into Section \ref{sectionatp}.  I am grateful to Michael Harris for helpful discussions, responses to my questions, and encouragement to complete this project, as well as enthusiastic suggestions for subsequent related projects.  

I would also like to thank the anonymous referee, for helpful comments and for suggesting that I add some remarks on functional equations (in analogue with the development in \cite{kaCM}).

I thank Mark Behrens for answering my questions about the conjectured implications of the $p$-adic Eisenstein series for homotopy theory.  I would also like to thank both Harris and Behrens for alerting me, in the first place, to the role of $p$-adic Eisenstein series in homotopy theory.

This project relies upon the ideas in Nicholas Katz's construction of Eisenstein measures, without which the current project would be impossible.  This paper also relies heavily upon the helpfully detailed presentation of $\ci$-Eisenstein series in Goro Shimura's papers; if Shimura's writing had not been so precise, I would have struggled much more with the non-$p$-adic portion of this project.

\section{Siegel Eisenstein series on certain unitary groups}\label{HoloESeriesSection}

\subsection{Unitary groups}\label{unitarygroups}
We now introduce the unitary groups with which we work throughout this paper.  The material in this section is a special case of the material in \cite[Section 21]{sh}, and the setup is also similar to that in \cite{EDiffOps} and \cite{SHL}.

Fix a CM field $K$, and denote its ring of integers by $\OK$.  Let $V$ be a vector space of dimension $n$ over $K$, and let $\langle v_1, v_2\rangle_V$ be a non-degenerate hermitian pairing on $V$.  
Let $-V$ denote the vector space $V$ with the hermitian pairing $-\langle v_1, v_2 \rangle_V$, and let
\begin{align*}
W & = 2V = V\oplus -V\\
\langle (v_1, v_2), (v_1', v_2')\rangle_W & = \langle v_1, v_1' \rangle_V+ \langle v_2, v_2' \rangle_{-V}\\
& = \langle v_1, v_1' \rangle_V- \langle v_2, v_2' \rangle_{V}.
\end{align*}
The hermitian pairing $\langle, \rangle_W$ defines an involution on $g\mapsto \tilde{g}$ on $\End_K(W)$ by 
\begin{align*}
\langle g(w), w'\rangle_W = \langle w, \tilde{g}(w')\rangle_W
\end{align*}
(where $w$ and $w'$ denote elements of $W$).  Note that this involution extends to an involution on $\End_{K\otimes_E R}\left(V\otimes_E R\right)$ for any $E$-algebra $R$.
We denote by $G$ the algebraic group such that for any $E$-algebra $R$, the $R$-points of $G$ are given by
\begin{align*}
G(R) = \left\{ g\in \gl_{K\otimes_E R}\left(W\otimes_E R\right)\middle|  g\tilde{g} = 1\right\}.
\end{align*}
Similarly, we define $H$ to be the algebraic group associated to $\langle, \rangle_V$ and $H'$ to be the algebraic group associated to $\langle, \rangle_{-V}$.  Note that $G(\IR)$ is of signature $(n,n)$.  Also, note that the canonical embedding
\begin{align*}
V\oplus V\hookrightarrow W
\end{align*}
induces an embedding
\begin{align*}
H(R)\times H'(R)\hookrightarrow G(R)
\end{align*}
for all $E$-algebras $R$.  When the $E$-algebra $R$ over which we are working is clear from context or does not matter, we shall write $U(W)$ for $G$, $U(V)$ for $H$, and $U(-V)$ for $H'$.

We write $W = V_d\oplus V^d,$ where $V_d$ and $V^d$ denote the maximal isotropic subspaces
\begin{align*}
V^d & = \left\{(v, v)|v\in V\right\}\\
V_d & = \left\{(v, -v)| v\in V\right\}.
\end{align*}
Let $P$ be the Siegel parabolic subgroup of $G$ stabilizing $V^d$ in $V_d\oplus V^d$ under the action of $G$ on the right.  Denote by $M$ the Levi subgroup of $P$ and by $N$ the unipotent radical of $P$.

\paragraph{Two convenient choices of bases}
Let $e_1, \ldots, e_n$ be an orthogonal basis for $V$ and $\phi$ be the matrix for $\langle, \rangle_V$ with respect to $e_1, \ldots, e_n$.  For $i= 1, \ldots, n$, let $e_{i+n} = e_i$.  Then the matrix for $\langle, \rangle_{2V}$ with respect to the basis $e_1, \ldots, e_{2n}$ is
\begin{align*}
\omega = \begin{pmatrix}\phi&0\\ 0& -\phi\end{pmatrix}.
\end{align*} With respect to the basis $e_1, \ldots, e_{2n}$, $(g, h)\in U(V)\times U(-V)$ is the element $\diag(g, h)\in U(2V)$.  Let $\alpha$ be a totally imaginary element of $K$, and let
\begin{align}
\eta = \begin{pmatrix}0& -1_n\\ 1_n&0\end{pmatrix}\label{etaweyl}\\
S = \begin{pmatrix}1_n & -\frac{\alpha}{2}\phi\\ -1_n & -\frac{\alpha}{2}\phi\end{pmatrix}.\nonumber
\end{align}
Note that
$S\eta { }^t\bar{S} = \alpha\omega.$
Using the change of basis given by $S$, we see that $U(W\otimes_E \IR)$ is isomorphic to 
\begin{align*}
G(\eta) = \left\{g\in\gl_{2n}\left(\IC\right)\middle| g\eta { }^t\bar{g} = \eta\right\}.  
\end{align*}
(Later, we shall refer to the adelic points of $G(\eta)$, by which we will mean
\begin{align*}
\left\{g\in\gl_{2n}\left(\adeles_K\right)\middle| g\eta { }^t\bar{g} = \eta\right\},
\end{align*}
where $\bar{g}$ denotes the element obtained by applying the involution coming from $K$ to $g$.)
In particular, letting $G(\omega) = \left\{g\in\gl_{2n}\left(\IC\right)\middle| g\omega { }^t\bar{g} = \omega\right\}$, we have
\begin{align*}
S^{-1} G(\omega) S = G(\eta),\\
\mbox{and}\\
(g, h)\mapsto S^{-1}\diag(g, h)S &= \begin{pmatrix} \frac{g+h}{2} & \frac{-\alpha(g-h)\phi}{4}\\ -\alpha^{-1}\phi^{-1}(g-h) & \frac{\phi^{-1}(g+h)\phi}{2}\end{pmatrix}\\
& = \begin{pmatrix}\frac{g+h}{2} & \frac{-\alpha(g-h)\phi}{4}\\ -\alpha^{-1}\phi^{-1}(g-h) & \frac{\overline{{ }^t (g+h)^{-1}}}{2}\end{pmatrix}\in G(\eta)
\end{align*}
gives an embedding of $G(\phi)\times G(-\phi)$ into $G(\eta)$.  Furthermore, $S^{-1}P(\omega)S = P$, where $P = P(\eta)$ is the parabolic  subgroup of $G(\eta)$ stabilizing $V^d$ in $V_d\oplus V^d$ under the action of $G$, and $P(\omega)$ is the stabilizer in $G(\omega)$ of $V^d = \left\{(v, v)\middle| v\in V\right\}$ inside of $V\oplus -V$ under action on the right.  The Eisenstein series in this paper are on $G(\eta)$ with respect to the parabolic subgroup $P = P(\eta)$ stabilizing $V^d$ (so they are Siegel Eisenstein series).

%
%
%
%

Note that the choice of basis for $V$ over $K$ fixes an isomorphism
\begin{align}\label{leviglv}
M\isomto \gl_K(V),
\end{align}
identifying $M$ with a subgroup of $G$ consisting of all elements of the form $\diag({ }^t\bar{h}^{-1}, h)$ with $h$ in $\gl_K(V)$.  (Also, note that the isomorphism $M\isomto \gl_K(V)$ extends to an isomorphism $M(R)\isomto\gl_{K\otimes_E R}(V\otimes R)$ for all $E$-algebras $R$.)
%
%

\paragraph{Choice of a Shimura datum}
Let $E$ be the maximal totally real subextension of the CM field $K$.  Our setup is as in \cite[Section 1.2]{SHL}; we review the details relevant to our construction.  We fix a Shimura datum $(GU(2V), X(2V))$, and a corresponding Shimura variety $Sh(W)$, according to the conditions in \cite{SHL} and in \cite{EDiffOps}.  Similarly, we attach Shimura data and Shimura varieties to $G(U(V)\times U(V))$.  Note that the symmetric domain $X(2V)$ is holomorphically isomorphic to the tube domain $\hn$ consisting of $[E: \IQ]$ copies of 
\begin{align*}
\left\{z\in M_{n\times n}(\IC)\middle| i({ }^t\bar{z}-z)>0\right\}.  
\end{align*}
Let $\mathcal{K}_{\infty}$ be the stabilizer in $G(\IR)$ of a fixed point 
\begin{align}\label{gimeldefn}
\gimel = (\gimel_{\sigma}) 
\end{align}
in $\prod_{\sigma: E\hookrightarrow \bar{\IQ}}\hn$, so we can identify $G(\IR)/K_{\infty}$ with $\hn$.  Given a compact open subgroup $\mathcal{K}$ of $G(\adeles_f)$, denote by $_{\mathcal{K}}Sh(W)$ the Shimura variety whose complex points are given by
\begin{align*}
G(\IQ)\backslash X\times G(\adeles_f)/\mathcal{K}.
\end{align*}
This Shimura variety is a moduli space for abelian varieties together with a polarization, an endomorphism, and a level structure (dependent upon the choice of $\mathcal{K}$).

\paragraph{Relationship between automorphic forms on $\hn$ and on the adelic points of $G(\eta)$}
%
%
%
Let $\mathcal{K} =\mathcal{K}_\infty\mathcal{K}_f$ be a compact open subgroup of the adelic points of $G(\eta)$, with $\mathcal{K}_\infty$ the stabilizer of $\gimel$. Let $\mathfrak{K}_{\infty}$ denote the complexification of $\mathcal{K}_\infty$ (so $\mathfrak{K}_\infty\subseteq \gln(\IC)\times\gln(\IC)$).  Let $\rho$ be a representation of $\mathfrak{K}_{\infty}\prod_{v\ndivides \infty}\mathcal{K}_v$.  Note that $_{\mathcal{K}}Sh(W)$ consists of copies of $\hn$.  Automorphic forms $f(g)$ of weight $\rho$ viewed as functions on the adelic points of $G(\eta)$ and automorphic forms $\tilde{f}(z)$ of weight $\rho$ viewed as functions on $\hn$ are related via
\begin{align*}
f\left(g\right) \leftrightarrow \rho\left(\left(\left(C\gimel + D\right), \left(\bar{C}{ }^t\gimel + \bar{D}\right)\right)k\right)\tilde{f}(g_{\infty}\gimel),
\end{align*}
for $g = \gamma g_{\infty}k$, with $g_{\infty} = \begin{pmatrix}A & B\\ C& D\end{pmatrix}\in G_{\infty}(\IR)$, $k\in\prod_{v\ndivides\infty}\mathcal{K}_v$, and $\gamma\in G(\IQ)$.  

As explained precisely in \cite[Section 22]{sh}, an automorphic form $f$ on $G(\eta)$ pulls back to an automorphic form $f(g, h)$ on $U(V)\times U(-V)$ which satisfies an automorphy property in terms of $g$ (independent of $h$) and an automorphy property in terms of $h$ (independent of $g$), which at each factor is completely determined by the automorphy property of $f$ on $G(\eta)$.  (The details of the spaces on which $G(V)$ acts are unnecessary for the current paper.  They are discussed in \cite[Sections 22, 6]{sh}.)
Later in the paper, we shall also work with automorphic forms (in particular, Eisenstein series) in a $p$-adic setting.


\subsection{Certain Eisenstein series}\label{eisensteinseriesprelims}
\paragraph{Choice of a CM type}
For the remainder of the paper, we fix a rational prime $p$ that is unramified in $K$ and such that each prime of the maximal totally real subfield $E\subseteq K$ dividing $p$ splits completely in the CM field $K$.  We also fix embeddings
\begin{align*}
\iota_\infty: &\bar{\IQ}\hookrightarrow\IC\\
\iota_p: &\bar{\IQ}\hookrightarrow\bar{\IQ}_p,
\end{align*}
and fix an isomorphism
\begin{align*}
\iota: \bar{\IQ}_p\isomto\IC
\end{align*}
satisfying $\iota\circ\iota_p = \iota_{\infty}$.  From here on, we identify $\bar{\IQ}$ with $\iota_p(\bar{\IQ})$ and $\iota_{\infty}(\bar{\IQ})$.

Fix a CM type $\Sigma$ for $K/\IQ$.  For each element $\sigma\in \Hom(E, \bar{\IQ})$, we also write $\sigma$ to denote the unique element of $\Sigma$ prolonging $\sigma:E\hookrightarrow\bar{\IQ}$ (when no confusion can arise).  For each element $x\in K$, denote by $\bar{x}$ the image of $x$ under the unique non-trivial element $\epsilon\in \Gal(K/E)$, and let $\bar{\sigma} = \sigma\circ\epsilon$.

Given an element $a$ of $E$, we identify it with an element of $E\otimes\IR$ via the embedding
\begin{align}\label{embconv}
E&\hookrightarrow E\otimes\IR\\
a&\mapsto (\sigma(a))_{\sigma\in\Sigma}.
\end{align}
We identify $a\in K$ with an element of $K\otimes \IC\isomto (E\otimes\IC)\times (E\otimes \IC)$ via the embedding
\begin{align*}
K&\hookrightarrow K\otimes\IC\\
a&\mapsto \left((\sigma(a))_{\sigma\in\Sigma}, (\bar{\sigma}(a))_{\sigma\in\Sigma}\right).
\end{align*}

\subsubsection{A unitary Hecke character of type $A_0$}\label{unitaryheckecharacter}
Let $\mathfrak {m}$ be an integral ideal in the totally real field $\Oe$ that divides $p^\infty$.  Let $\chi$ be a unitary Hecke character of type $A_0$
\begin{align*}
\chi: \adeles_K^\times\rightarrow \IC^\times
\end{align*}
of conductor $\mathfrak{m}$, i.e. 
\begin{align*}
\chi_v(a) = 1
\end{align*}
for all finite primes $v$ in $K$ and all $a\in K_v^\times$ such that
\begin{align*}
a\in 1+\mathfrak{m}_v{\OK}_v.
\end{align*}
Let $\nu(\sigma)$ and $k(\sigma)$, $\sigma\in\Sigma$, denote integers such that the infinity type of $\chi$ is given by
\begin{align}\label{infinitytypeknu}
\prod_{\sigma\in\Sigma}\sigma(b)^{-k(\sigma)-2\nu(\sigma)}\left(\sigma(b)\bar{\sigma}(b)\right)^{\frac{k(\sigma)}{2}+\nu(\sigma)}.
\end{align}
\begin{rmk}
Given an element $a\in K$, we associate $a$ with an element of $K\otimes \IR$, via the embedding
\begin{align*}
a\mapsto (\sigma(a))_{\sigma\in \Sigma},
\end{align*}
so for $a\in K$, 
\begin{align*}
\left(\prod_{v\divides\infty}\chi_v\right)(a) = \prod_{\sigma\in\Sigma}\chi_v(\sigma_v(a)).
\end{align*}
\end{rmk}
\paragraph{Conventions for adelic norms}
Let $\left|\cdot\right|_E$ denote the adelic norm on $E^\times\backslash\adeles_E^\times$ such that for all $a\in \adeles_E^\times$, 
\begin{align*}
\left|a\right|_E = \prod_v\left|a\right|_v,
\end{align*}
where the product is over all places of $E$ and where the absolute values are normalized so that
\begin{align*}
\left|v\right|_v & = q^{-1},\\
q & = \mbox{ the cardinality of } {\Oe}_v/v{\Oe}_v,
\end{align*}
for all non-archimedean primes $v$ of the totally real field $E$.  Consequently for all $a\in E$,
\begin{align*}
\prod_{v|\infty}|a|_v^{-1} = \prod_{v\in\Sigma}\sigma_v(a)\Sgn(\sigma_v(a)),
\end{align*}
where the product is over all archimedean places $v$ of the totally real field $E$.  We denote by $\left|\cdot\right|_K$ the adelic norm on $K^\times\backslash\adeles_K^\times$ such that for all $a\in \adeles_K^\times$,
\begin{align*}
\left|a\right|_K = \left|a\bar{a}\right|_E.
\end{align*}
For $a\in K$ and $v$ a place of $E$, we let
\begin{align*}
\left|a\right|_v = \left|a\bar{a}\right|_v^{\frac{1}{2}}.
\end{align*}

\subsubsection{Siegel Eisenstein series}\label{paraboliccharacter}
Let $\chi$ be a unitary Hecke character that meets the conditions of Section \ref{unitaryheckecharacter}.  For any $s\in \IC$, we view $\chi\cdot\left|\cdot\right|_K^{-s}$ as a character of the parabolic subgroup $P(\adeles_E) = M(\adeles_E)N(\adeles_E)$ via the composition of maps
\begin{tiny}
\begin{align*}
P\left(\adeles_E\right)\xrightarrow{\mod N\left(\adeles_E\right)} M\left(\adeles_E\right)&\xrightarrow{\mbox{the iso in \eqref{leviglv}}} \gl_{\adeles_K}(V\otimes_E\adeles_E)
\xrightarrow{\det} \adeles_K^\times \xrightarrow{\chi\left|\cdot\right|_K^{-s}} \IC^\times.
\end{align*}
\end{tiny}
Consider the induced representation
\begin{align}\label{inducedrep}
I(\chi, s) = \Ind_{P(\adeles_E)}^{G(\adeles_E)} (\chi\cdot\left|\cdot\right|_K^{-s})\isomto \otimes_v\Ind_{P(E_v)}^{G(E_v)}\left(\chi_v\left|\cdot\right|_v^{-2s}\right),
\end{align}where the product is over all places of $E$.


\begin{rmk}Although it is in some contexts conventional to use normalized induction and consider $\chi\cdot\left|\cdot\right|_K^{s-{n/2}}$ instead of $\chi\cdot\left|\cdot\right|_K^s$, we leave out the factor $\left|\cdot\right|_K^{-{n/2}}$ for now.  (It just comes from the modulus character and will reappear when we compute the Fourier coefficients.) The reader who prefers to include $n/2$ in this expression may just regard it as absorbed into $s$.  
\end{rmk}

Given a section $f\in I(\chi, s)$, the Siegel Eisenstein series associated to $f$ is the $\IC$-valued function of $G$ defined by
\begin{align*}
E_f(g) = \sum_{\gamma\in P(E)\backslash G(E)}f(\gamma g)
\end{align*}
This function converges for $\Re(s)>0$ and can be continued meromorphically to the entire plane.
\begin{rmk}
If we were working with normalized induction, then the function would converge for $\Re(s)>\frac{n}{2}$, but we have absorbed the exponent $\frac{n}{2}$ into the exponent $s$. (Our choice not to include the modulus character at this point is equivalent to shifting the plane on which the function converges by $\frac{n}{2}$.)
\end{rmk}
All the poles of $E_f$ are simple and there are at most finitely many of them.  Details about the poles are given in \cite{tan}.

\subsubsection{Conditions imposed on the choice of a Siegel section}\label{conditions}
We continue to work with a unitary Hecke character $\chi$ meeting the conditions of Section \ref{unitaryheckecharacter}, which has infinity type given by \eqref{infinitytypeknu}.  Let $k(\sigma)$ and $\nu(\sigma)$ be as in \eqref{infinitytypeknu}.  

In order to construct an Eisenstein series that behaves nicely with respect to $p$-adic interpolation, we shall work with a section $f\in I(\chi, s)$ that meets the following conditions (among other conditions that we will specify as they become relevant):
\begin{condition}\label{condA}
There is a compact open subgroup $U$ of $G(\adeles_E)$ and a $\bar{\IQ}^\times$-valued character $\rho$ of $\prod_{v\divides p}U_v$ such that 
for all $u\in U$ and $h\in G(\adeles_E)$,
\begin{align*}
E_f(hu) & = \rho((u_v)_{v\divides p})J_{(u_v){v\divides\infty}}^{k, \nu}(\gimel)^{-1}E_f(h),
\end{align*}
where $k = (k(\sigma))_{\sigma\in\Sigma}$, $\nu = (\nu(\sigma))_{\sigma\in\Sigma}$, and  $J_{\alpha}^{k, \nu}(z)$ is defined for all $z\in \hn$ and 
\begin{align*}
\alpha &= \begin{pmatrix}a& b\\ c&d\end{pmatrix} = \left(\begin{pmatrix}a_{\sigma}& b_{\sigma}\\ c_{\sigma}& d_{\sigma}\end{pmatrix}\right)_{\sigma\in\Sigma}\in \prod_{\sigma\divides\infty}G(\IR)
\end{align*}
by
\begin{align}
J_{\alpha}^{k, \nu}(z) &: = \prod_{\sigma\in\Sigma}(\det\alpha_{\sigma})^{\nu(\sigma)}\det(c_{\sigma}z_{\sigma}+d_{\sigma})^{k(\sigma)}\label{jknudef2}\\
& = \prod_{\sigma\in\Sigma}\det(c_{\sigma}z_{\sigma}+d_{\sigma})^{k(\sigma)+\nu(\sigma)}\prod_{\sigma\in\Sigma}\det\left(\bar{c}_{\sigma}{ }^tz_{\sigma}+\bar{d}_{\sigma}\right)^{-\nu(\sigma)}\nonumber.
\end{align}
\end{condition}

\begin{condition}\label{condC}The section $f$ factors as $f = \otimes_vf_v$ with $f_v\in\Ind_{P(E_v)}^{G(E_v)}\left(\chi_v\left|\cdot\right|_v^{2s}\right)$ for each place $v$ of $E$.  (This condition allows us to express the Fourier coefficients as a product of local factors, which is important for our approach to $p$-adic interpolation of the the Fourier coefficients.)
\end{condition}

\begin{rmk}
Many different choices for the local sections $f_v$ lend themselves nicely to construction of a $p$-adic Eisenstein measure.  The most constrained choices are at $p$ and at $\infty.$  
\end{rmk}

\subsubsection{Preliminaries on Fourier expansions}\label{consequencesAB}

\paragraph{Exponential characters}
Before explaining the consequences of Conditions \ref{condA} and \ref{condC}, we need to introduce some conventions for exponential characters.  For each archimedean place $v\in\Sigma$, denote by $\e_v$ the character of $E_v$ (i.e $\IR$) defined by
\begin{align*}
\e_v(x_v) = e^{(2\pi i x_v)}
\end{align*}
for all $x_v$ in $E_v$.  Denote by $\e_{\infty}$ the character of $E\otimes \IR$ defined by
\begin{align*}
\e_{\infty}((x_v)_{v\in \Sigma}) = \prod_{v\divides\infty}\e_v(x_v).
\end{align*}
Following our convention from \eqref{embconv}, we put
\begin{align*}
\e_{\infty}(a) = \e_{\infty}((\sigma(a))_{\sigma\in\Sigma}) = \e^{2\pi i\tr_{E/\IQ}(a)}
\end{align*}
for all $a\in E$.
For each finite place $v$ of $E$ dividing a prime $q$ of $\ZZ$, denote by $\e_v$ the character of $E_v$ defined for each $x_v\in E_v$ by
\begin{align*}
\e_v(x_v) = e^{-2\pi i y}
\end{align*}
where $y$ is an element of $\IQ$ such that $\tr_{E_v/\IQ_q}(x_v)-y\in\ZZ_p$.  We denote by $\e_{\adeles_E}$ the character of $\adeles_E$ defined by
\begin{align*}
\e_{\adeles_E}(x) = \prod_v\e_v\left(x_v\right)
\end{align*}
for all $x = \left(x_v\right)\in \adeles_E$.  
\begin{rmk}
Note that for a $a\in E$, we identify $a$ with the element $(\sigma_v(a))_v\in\adeles_E$, where $\sigma_v: E\hookrightarrow E_v$ is the embedding corresponding to $v$.  Following this convention, we put
\begin{align}\label{ewarning}
\e_{\adeles_E}(a) = \prod_v\e_v(\sigma_v(a)).
\end{align}
for all $a\in E$.
\end{rmk}

\paragraph{Consequences of Conditions \ref{condA} and \ref{condC}}
For any subring $R$ of $K\otimes_E E_v$, with $v$ a place of $E$, let $\hern(R)$ denote the space of $n\times n$-matrices with entries in $R$.  

By \cite[Proposition 18.3]{sh}, Condition \ref{condA} guarantees that $E_f$ has a Fourier expansion such that for all $h\in \gln(K)$ and $m\in\hern(K)$
\begin{align*}
E_f\left(\begin{pmatrix}1 & m\\ 0& 1\end{pmatrix}\begin{pmatrix}{ }^t\bar{h}^{-1} & 0\\ 0 & h\end{pmatrix}\right) = \sum_{\beta\in\hern(K)}c(\beta, h; f)\e_{\adeles_E}\left(\tr\left(\beta m\right)\right),
\end{align*}
with $c(\beta, h; f)$ a complex number dependent only on the choice of section $f$, the hermitian matrix $\beta\in \hern(K)$, $h_v$ for finite places $v$, and $\left(h\cdot  { }^t\bar{h}\right)_v$ for archimedean places $v$ of $E$.

By \cite[Sections 18.9, 18.10]{sh} and Condition \ref{condC}, the Fourier coefficient is a product of local Fourier coefficients determined by the local sections $f_v$.  More precisely, for each $\beta\in\hern(K)$,
\begin{small}
\begin{align}
c(\beta, h; f) &= C(n, K)\prod_vc_v(\beta_v, h_v; f)\nonumber\\
c_v(\beta, h; f) &= \nonumber\\
\prod_v&\int_{\hern(K\otimes F_v)} f_v\left(\begin{pmatrix}0& -1\\ 1& 0\end{pmatrix}\begin{pmatrix}1 & m_v\\ 0& 1\end{pmatrix}\begin{pmatrix}{ }^t\bar{h_v}^{-1} & 0\\ 0 & h_v\end{pmatrix}\right)\e_v(-\tr(\beta_v m_v))dm_v\nonumber\\
C(n, K)  &= 2^{n(n-1)[E:\IQ]/2}\left|D_F\right|^{-n/2}\left|D_K\right|^{-n(n-1)/4}\label{CnK},
\end{align}
\end{small}
where $D_F$ and $D_K$ are the discriminants of $K$ and $F$ respectively, $\beta_v = \sigma_v(\beta)$ for each place $v$ of $F$, and $d_v$ denotes the Haar measure on $\hern(K_v)$ such that:
\begin{align}
&\int_{\hern\left(\OK\otimes_F F_v\right)} d_vx = 1, \mbox{ for each finite place $v$ of F}\nonumber\\
d_vx&: = \left|\bigwedge_{j=1}^{n}dx_{jj}\bigwedge_{j<k}\left(2^{-1}dx_{jk}\wedge d\bar{x}_{jk}\right)\right|, \mbox{ for each archimedean place $v$ of $E$}.\label{matrixnote}
\end{align}
(In Equation \eqref{matrixnote}, $x$ denotes the matrix whose $ij$-th entry is $x_{ij}$.)

\begin{warning}\label{betawarning}
Given $\beta\in \hern(K)$, we shall sometimes write $\e_v(\beta)$ when we mean $\e_v(\sigma_v(\beta))$, in keeping with the convention of Equation \eqref{ewarning}.  (This agrees with Shimura's convention in \cite{sh} and \cite{shconfluent}.)  Note that this notational ambiguity in Shimura's papers seems to have led to an error in (3.3.5.3) and in computation of the Fourier coefficient at $\infty$ in \cite{SHL}.  In particular, the factor $\det\beta^{[E:\IQ]}$ in \cite[Equation (3.3.5.3)]{SHL} should be $\det(\prod_{\sigma\in\Sigma}\sigma(\beta))$.  This correction is necessary for our construction of the Eisenstein measure.
\end{warning}

\paragraph{Fourier coefficients at points on the Levi subgroup}
We continue to work with a section $f=\otimes_v f_v$ that lies in the induced representation \eqref{inducedrep} and meets all the above conditions.
\begin{lem}\label{mavsm1}
For each $h\in \gln(\adeles_K)$ and $\beta\in\hern(K)$,
\begin{align}\label{fcoeffsat1}
c(\beta, h; f) = \chi({ }^t\bar{h}^{-1})\left|\det\left({ }^t\bar{h}^{-1}\cdot h^{-1}\right)\right|_E^{n-s} c({ }^t\bar{h}^{-1}\beta h^{-1}, 1_n; f).
\end{align}
\end{lem}
\begin{proof}
As in Equation \eqref{etaweyl}, let $\eta = \begin{pmatrix}0 & -1_n\\ 1_n & 0\end{pmatrix}.$
Let $m(h)$ denote the matrix $\begin{pmatrix}{ }^t\bar{h}^{-1} & 0\\ 0 & h\end{pmatrix}$.  Observe that for any $n\times n$ matrix $m$,
\begin{align*}
\eta \cdot m(h)\cdot \eta^{-1} & = m({ }^t\bar{h}^{-1})\\
m(h)^{-1}\cdot \begin{pmatrix}1&m\\ 0 &1\end{pmatrix}\cdot m(h) &= \begin{pmatrix}1 & { }^t\bar{h} m h\\ 0 & 1\end{pmatrix}.
\end{align*}
Therefore,
\begin{align*}
\eta \cdot \begin{pmatrix}1 & m \\ 0 & 1\end{pmatrix}\cdot m(h) & = \left(\eta \cdot m(h)\cdot \eta^{-1}\right)\eta \left(m(h)^{-1}\begin{pmatrix}1 & m\\ 0 & 1\end{pmatrix} m(h)\right)\\
 & = m({ }^t\bar{h}^{-1}) \eta \begin{pmatrix}1 & { }^t\bar{h} m h\\ 0 & 1\end{pmatrix}.
\end{align*}
So for any place $v$ of $E$ and section $f_v\in \Ind_{P(E_v)}^{G(E_v)}(\chi, s)$,
\begin{align}\label{leviunimp}
f_v\left(\eta\begin{pmatrix}1 & m \\ 0 & 1\end{pmatrix}m(h_v)\right) = \chi_v\left({ }^t\bar{h_v}^{-1}\right)\left|{ }^t\bar{h_v}^{-1}\right|^{-2s}_vf_v\left(\eta \begin{pmatrix}1 & { }^t\bar{h}_v m h_v\\ 0 & 1\end{pmatrix}\right) \end{align}
The lemma now follows from Equation \eqref{leviunimp} and the fact that the Haar measure $d_v$ satisfies $d_v(h_vx{ }^t\bar{h_v}) = \left|\det\left({ }^t\bar{h_v}\cdot h_v\right)\right|_v^{n}d_v(x)$ for each place $v$ of $E$. 
\end{proof}

\subsubsection{The sections at archimedean places of $E$}

As in Condition \ref{condA}, let $k = (k(\sigma))_{\sigma\in\Sigma}$ and $\nu = (\nu(\sigma))_{\sigma\in\Sigma}$, and let $J_{\alpha}^{k, \nu}$ be defined as in Equation \eqref{jknudef2}.  Given a unitary Hecke character $\chi$ meeting the conditions of Section \ref{unitaryheckecharacter}, a complex number $s$, and $\gimel$ as in Equation \eqref{gimeldefn},  we will always take the section at the archimedean places to be of the form
\begin{align*}
f_{\infty}^{k, \nu}(\bullet; \gimel, \chi, s)\in\otimes_{v\divides\infty}\Ind_{P(E_v)}^{G(E_v)}\left(\chi_v\left|\right|_v^{-2s}\right),
\end{align*}
where
\begin{align*}
f_{\infty}^{k, \nu}(\alpha;  \gimel, \chi, s) & =J_{\alpha}^{k, \nu}(\gimel)^{-1}(J_{\alpha}(\gimel)\overline{J_{\alpha}(\gimel)})^{\frac{k}{2}-s},\\
\left(J_{\alpha}(z)\overline{J_{\alpha}(z)} \right)^{\frac{k}{2}-s}& : = \prod_{\sigma\in\Sigma}\left(\det(c_{\sigma} z_{\sigma}+ d_{\sigma})\overline{\det(c_{\sigma} z_{\sigma}+ d_{\sigma})}\right)^{\frac{k(\sigma)}{2}-s},
\end{align*}
for each $\alpha = \begin{pmatrix}a_{\sigma} & b_{\sigma}\\ c_{\sigma} & d_{\sigma}\end{pmatrix} \in \prod_{v\divides\infty}G(E_v)$.  The relationship between $f^{k,\nu}_{\infty}(\alpha; \gimel, \chi, s)$ and $f^{k, \nu}_{\infty}(\alpha; i\cdot 1_n, \chi, s)$ is given by
\begin{align}\label{gimelvsi}
f^{k, \nu}_{\infty}(\alpha; \gimel, \chi, s) = f^{k, \nu}_{\infty}(\alpha g^{-1}; i 1_n, \chi, s)f^{k, \nu}_{\infty}(g^{-1}; i 1_n, \chi, s)^{-1},
\end{align}
where $g$ is any element of $G$ such that $g\gimel = i 1_n$.

\subsubsection{The Fourier coefficients at archimedean places of $E$}\label{sectionsinfty}

Recall that by Equation \eqref{fcoeffsat1}, the Fourier coefficients $c(\beta, h; f)$ are completely determined by their values at $h= 1_n$.  Following Shimura, we give the Fourier coefficients in terms of the section $f^{k, \nu}_{\infty}(\bullet; i 1_n, \chi, s)$ and note that Equation \eqref{gimelvsi} relates the Fourier coefficients in terms of $f^{k, \nu}_{\infty}(\bullet;\gimel , \chi, s) = \otimes_v f_v^{k, \nu}$ to those in terms of $f^{k, \nu}_{\infty}(\bullet; i 1_n, \chi, s)$.
The non-zero Fourier coefficients of the Eisenstein series at the archimedean primes, evaluated at $m(1_n)$ are given by the following product over all archimedean places of $E$:
\begin{small}
\begin{align*}
\prod_{v\divides\infty}c_v(\beta,& 1_n; f_v^{k, \nu}\left(\bullet; i1_n, \chi, s\right))\\ = \prod_{v\divides\infty}&\int_{\hern(\IC)}f_v^{k, \nu}\left(\begin{pmatrix}0& -1\\ 1& 0\end{pmatrix}\begin{pmatrix}1 & m_v \\ 0 & 1\end{pmatrix}\right) \e{\left( -\tr\left(\sigma_v(\beta) m_v\right)\right)}dm_v\\
 = \prod_{v\divides\infty}&\int_{\hern(\IC)}f_v^{k, \nu}\left(\begin{pmatrix}0& -1\\ 1 & m_v\end{pmatrix}\right) \e{\left( -\tr\left(\sigma_v(\beta) m_v\right)\right)}dm_v\\
 = \prod_{v\divides\infty}&\int_{\hern(\IC)}\det (i + m_v)^{-\frac{k(\sigma_v)}{2}-s}\det\left(-i+m_v\right)^{\frac{k(\sigma_v)}{2}-s}\e(\tr(-\sigma_v(\beta) m_v))dm_v.
\end{align*}
\end{small}
By \cite[Lemma 18.12]{sh}, which is proved in \cite{shconfluent}, each of the above local integrals converges at least for $\Re\left(\frac{k(\sigma)}{2}+s\right)>2n-1$ and can be meromorphically continued to a meromorphic function of $(k(\sigma), s)\in \IC\times\IC$.
When there is an integer $k$ such that 
\begin{align*}
s = \frac{k}{2} = \frac{k(\sigma)}{2}  \mbox{ for all $\sigma\in\Sigma$}
\end{align*}
(i.e. when $J_{\alpha}^{k, \nu}(z)^{-1}(J_{\alpha}(z)\overline{J_{\alpha}(z)})^{\frac{k}{2}-s}$ is a holomorphic function of $z\in \hn$), \cite[Equation (7.12)]{sheseries} describes the archimedean Fourier coefficients precisely:
\begin{small}
\begin{align}
c_v&\left(\beta, 1_n; f_v^{k, \nu}\left(\bullet; i1_n, \chi, \frac{k}{2}\right)\right)\nonumber \\
 &= 2^{(1-n)n}i^{-nk}(2\pi)^{nk}\left(\pi^{n(n-1)/2}\prod_{t=0}^{n-1}\Gamma(k-t)\right)^{-1}\sigma_v(\det\beta)^{k-n}\e\left(i\tr(\sigma_v(\beta))\right)\label{coeffinfty},
\end{align}
\end{small}
for each archimedean place $v$ of $E$.  Observe that when $k>n$, 
\begin{align*}
\prod_{v\divides \infty}c_v\left(\beta, h; f_v^{k, \nu}\left(\bullet; i1_n, \chi, \frac{k}{2}\right)\right) = 0,
\end{align*}
unless $\det(\beta)\neq 0$ and $\det(h)\neq 0$, i.e. unless $\beta$ is of rank $n$.

\subsubsection{Conventions at places of $E$ dividing $p$}
We now establish certain conventions that we will use in Section \ref{sectionatp}, when we define the Siegel sections at places of $E$ dividing $p$.

\paragraph{CM Types}
As explained in \cite[Section (5.1.10)]{kaCM}, the choice of the CM type $\Sigma$ is equivalent to the choice of a prime in $K$ over each prime in $E$ dividing $p$.  We write ${v\in\Sigma}$ to mean $v$ is a prime in $K$ dividing $p$ that induces a $p$-adic embedding $\iota_p\circ\sigma$ for some embedding $\sigma\in\Sigma$, and we associate $v\in\Sigma$ with this embedding $\sigma\in\Sigma$.  When no confusion can arise, we write $v$ to denote the prime in $\Sigma$ lying above a prime $v$ in $E$ dividing $p$.  For each element $x\in K$, denote by $\bar{x}$ the image of $x$ under the unique non-trivial element of $\epsilon\in\Gal(K/E)$, let $\bar{v}=\epsilon(v)$, and let $\bar{\sigma} = \sigma\circ\epsilon$.  Let 
\begin{align}\label{phisigma}
\phi_{\Sigma}: K\otimes\IQ_p\rightarrow E\otimes\IQ_p
\end{align}
be the projection obtained through the composition of maps
\begin{tiny}
\begin{align*}
K\otimes\IQ_p\twoheadrightarrow \prod_{v\in\Sigma}(v\mathrm{-adic}\mbox{ completion of }K)\isomto \prod_{v\in\Sigma}(v\mathrm{-adic}\mbox{ completion of }E) = E\otimes\IQ_p
\end{align*}
\end{tiny}
From here on, we identify
$K\otimes\IQ_p$ with $\left(E\otimes\IQ_p\right)\times\left(E\otimes\IQ_p\right)$ via the canonical ring isomorphism
\begin{align}\label{okoeid}
K\otimes\IQ_p&\isomto\left(E\otimes\IQ_p\right)\times\left(E\otimes\IQ_p\right)\\
x&\mapsto (\phi_{\Sigma}(x), \phi_{\Sigma}(\bar{x})) = ((\sigma(x))_{\sigma\in\Sigma}, (\sigma(\bar{x})_{\sigma\in\Sigma})\nonumber
\end{align}
\paragraph{Conventions for certain $p$-adic maps}\label{chiconvention}
The isomorphism \eqref{okoeid} induces an identification of $V\otimes\IQ_p$ with a module over $\left(E\otimes\IQ_p\right)\times\left(E\otimes\IQ_p\right)$.  Similarly, the isomorphism \eqref{okoeid} identifies modules over $\OK\otimes\ZZ_p$ with modules over $(\Oe\otimes\ZZ_p)\times(\Oe\otimes\ZZ_p)$.
Given continuous maps (most often characters) 
\begin{align*}
\chi_1, \chi_2: \gln (E\otimes\IQ_p)\rightarrow\IC_p^\times,
\end{align*}
we write $(\chi_1, \chi_2^{-1})$ to mean the function on $\gln(K\otimes\IQ_p)$ defined by
\begin{align}\label{chi1chi2inv}
a\mapsto \chi_1(\phi_{\Sigma}(a))\chi_2(\phi_{\Sigma}({ }^t\bar{a}^{-1})),
\end{align}
so
\begin{align*}
{ }^t\bar{a}^{-1} \mapsto \chi_1(\phi_{\Sigma}({ }^t\bar{a}^{-1}))\chi_2(\phi_{\Sigma}(a)).
\end{align*}
\subparagraph{Hermitian matrices}
In particular, if $\beta\in \gln\left(\OK\otimes\ZZ_p\right)$ is Hermitian, then
\begin{align*}
\beta\mapsto \chi_1(\phi_{\Sigma}(\beta))\chi_2(\phi_{\Sigma}(\beta^{-1})).
\end{align*}
We identify the space of Hermitian matrices with its image  under the map induced by $\phi_{\Sigma}$ (i.e project onto the first factor in the map \eqref{okoeid}).

\subparagraph{The unitary groups $U(V)$ and $U(W)$ at $p$}Let $H = U(V)$ and $G = U(W) =U(2V)$ (as in Section \ref{unitarygroups}).  As in Section \ref{unitarygroups}, fix a basis $e_1, \ldots, e_n$ for $V$ over $K$.  We identify $V^d$ and $V_d$ with $V$ via the isomorphisms
\begin{align*}
V^d\isomto &V\isomto V_d\\
\left(e_i, e_i\right)\mapsto &e_i\mapsto \left(e_i, -e_i\right)
\end{align*}
induced by the choice of basis $e_1, \ldots, e_n$.  Via the isomorphism \eqref{okoeid}, we identify $V\otimes\IQ_p$ with the direct product of two $n$-dimensional vector spaces over $E\otimes\IQ_p$, each of which is fixed under the action of the unitary group $H(E\otimes\IQ_p)$ and which are swapped under the involution induced by complex conjugation $\epsilon$ on $K$.  Similarly, we identify $W\otimes\IQ_p$ with the direct product of two $2n$-dimensional vector spaces over $E\otimes\IQ_p$, each of which is fixed under the action of the unitary group $G(E\otimes\IQ_p)$ and which are swapped under the involution induced by $\epsilon$.

For each prime $v\in\Sigma$, we use the following identifications, induced by the projection map $\phi_{\Sigma}$ defined in \eqref{phisigma}:
\begin{align*}
H(E_v)\isomto \gl_{K_v}(V\otimes K_v)\\
G(E_v)\isomto \gl_{K_v}(W\otimes K_v).
\end{align*}
For each prime $v\in\Sigma$, the choice of basis $e_1, \ldots, e_n$ together with the above isomorphism also induces identifications
\begin{align}
H(E_v)\isomto &\gln(K_v)\isomto \gln(E_v)\nonumber\\
G(E_v)\isomto &\gl_{2n}(K_v)\isomto \gl_{2n}(E_v)\label{ggl2n}.
\end{align}

\subparagraph{The parabolic and Levi subgroups at $p$} Let $v$ be a prime in $\Sigma$.  The image of $M(E_v)$ in $GL_{2n}(E_v)\isomto G(E_v)$ under the above identifications is
\begin{tiny}
\begin{align*}
\left\{\diag(B, A)\in \gl_{2n}(E_v)\middle| B, A\in \gln(E_v)\isomto \gl_{K_v}(V\otimes_K K_v)\right\}\isomto \gln(E_v)\times \gln(E_v).
\end{align*}
\end{tiny}
The isomorphism \eqref{leviglv} together with the isomorphism \eqref{okoeid}  gives identifications
\begin{align}\label{leviab}
M(E_v)&\isomto \gl_{K_v}(V\otimes_K K_v)\times \gl_{K_{\bar{v}}}(V\otimes_K K_{\bar{v}})\nonumber\\
&\isomto \gl_{E_v}(V\otimes E_v)\times \gl_{E_v}(V\otimes_E E_{v}).
\end{align}
Let $m=\diag(A, B)\in \gl_{2n}(E_v)$ be an element of the Levi subgroup $M(E_v)$.   The image of $m$ under the isomorphism \eqref{leviab} is $(A, { }^t B^{-1})$.  

The image of the parabolic subgroup $P(E\otimes\IQ_p)$ in $\gl_{2n}(E_v)$ under the identification \eqref{ggl2n} of $G(E\otimes\IQ_p)$ with $\gl_{2n}(E_v)$ is
\begin{align*}
\left\{\begin{pmatrix}B & C\\ 0 & A\end{pmatrix}\middle| A, B\in \gl_n(E_v), C\in M_{n\times n}(E_v) \right\}.
\end{align*}
Let $\alpha$ be an element of $P(E\otimes\IQ_p)$ whose image in $\gl_{2n}(E\otimes\IQ_p)$ under the above identification is $\begin{pmatrix}B & C\\ 0 & A\end{pmatrix}$.  Let $\chi_p$ denote the restriction of the Hecke character $\chi$ to
\begin{align*}
K\otimes \IQ_p \stackrel{\eqref{okoeid}}{=} \prod_{v\in\Sigma}K_v\times \prod_{v\in\Sigma} K_{\bar{v}} = \prod_{v\in\Sigma}E_v\times \prod_{v\in\Sigma} E_v.
\end{align*}
Writing
\begin{align*}
\chi_p = (\chi_1, \chi_2^{-1})\\
\left|\cdot \right|^{-s}_p = \prod_{v\divides p}\left|\cdot \right|^{-s}_v
\end{align*}
as in \eqref{chi1chi2inv} and evaluating $\chi_p\cdot \left|\cdot \right|_p^{-s}$ on $\alpha$ as in Section \ref{paraboliccharacter}, we have that 
\begin{align*}
\chi_p(\alpha)\left|\alpha\right|_p^{-s} = \chi_1(\det A)\chi_2(\det B)\prod_{v\in\Oe, v\divides p}\left|\det (B^{-1} A)\right|_v^{-s}.
\end{align*}

%
\subsubsection{Siegel sections at $p$}\label{sectionatp}

For the reasons explained in Item \eqref{padiccorrections} of Section \ref{CorrectionsClarifications}, we modify the Siegel section in \cite{SHL}.  Since different sections at $p$ are likely best in different contexts, we also give a recipe for constructing certain kinds of sections at $p$ that can be used to produce an Eisenstein measure in Section \ref{EisensteinMeasureSection}.  At the end of this section, we discuss some natural generalizations to our construction.  


Given the character $\chi_p:= \prod_{v|p}\chi_v$ of $(K\otimes\IQ_p)^\times\cong (E\otimes\IQ_p)^\times\times(E\otimes\IQ_p)^\times$ (where the isomorphism is as in \eqref{okoeid}), we write $\chi_p= (\chi_1, \chi_2^{-1})$, using the convention established in Section \ref{chiconvention}.

Let $v\subseteq K$ be a prime in $\Sigma$.  For any $K$-vector space $U$, denote by $U_v$ the vector space $U\otimes_K K_v$.  To each Schwartz function
\begin{align*}
\Phi_v: \Hom_{K_v}(V_v, V_{v, d}\oplus V_v^d)\rightarrow K_v, 
\end{align*} we attach a section
\begin{align*}
f^{\Phi_v}\in \Ind_{P(E_{v})}^{G(E_{v})}\left(\chi_{v}\cdot\left|\cdot\right|_{v}^{-2s}\right) 
\end{align*}
as follows.
Consider the decomposition
\begin{align*}
\Hom_{K_v}(V_v, W_v) &= \Hom_{K_v}(V_v, V_{d, v})\oplus\Hom_{K_v}(V_v, V^d_v)\\
X&= (X_1, X_2).
\end{align*}
Let
\begin{align*}
\mathbf{X} &= \left\{X\in\Hom_{K_v}(V_v, W_v)|X(V_v) = V_v^d\right\}\\
&= \left\{(0, X)|X: V_v\rightarrow V^d_v \mbox{ is an isomorphism}\right\}.
\end{align*}
Note that for any $X\in \mathbf{X}$, the composition $V_v\xrightarrow{X} V_v^d\isomto V_v$ is an isomorphism of $V$ with itself.  So we may identify $\mathbf{X}$ with $\Aut_{K_v}(V_v)$.  Via the conventions established in Section \ref{chiconvention}, we identify each $g\in G(E_v)$ with its image in $\gl_{K_v}(W_v)\isomto GL_{2n}(K_v)\rightarrow GL_{2n}(E_v)$.
We define the section $f^{\Phi_v}\in \Ind_{P(E_{v})}^{G(E_{v})}\left(\chi_{v}\cdot\left|\cdot\right|_{v}^{-2s}\right)$ by\footnote{A typo in the exponents in Equation \eqref{sectionininducedrepnatp} appeared in the first three versions of this paper on the arxiv.  The correct expression appears here.}
\begin{align}\label{sectionininducedrepnatp}
f^{\Phi_v}(g): = \chi_{2, v}(\det g)\left|\det g\right|_{v}^{2s}\int_{\mathbf{X}}\Phi(Xg)\chi_{1, v}^{-1}\chi_{2, v}(\det X)\left|\det X\right|_{v}^{4s}d^{\times}X.
\end{align}
It is a simple computation to check that, indeed,
$f^{\Phi_v}\in \Ind_{P(E_v)}^{G(E_v)}\left(\chi_v\left|\cdot \right|_v^{-2s}\right)$.

\paragraph{Fourier expansion at $p$}\label{FExpnatp}
Our choice of a section well-suited to $p$-adic interpolation relies upon Lemma \ref{pfflemma}, which we state below.

Given a Schwartz function
\begin{align*}F: \Hom_{K_v}(V_v, V_{d, v})\oplus\Hom_{K_v}(V_v, V^d_v)&\rightarrow \IC\\
(X_1, X_2)&\mapsto F(X_1, X_2),
\end{align*}
we denote by $PF(X_1, X_2)$ the Fourier transform of $F$ in the second variable, i.e.
\begin{align*}
PF(X_1, X_2) = \int_{M_n(E_v)}F(X_1, Y)\e_p^{\tr(-Y\cdot { }^tX_2)}dY.
\end{align*}
Note that
\begin{align*}
PPF(X_1, X_2) = F(X_1, -X_2).
\end{align*}

\begin{lem}\label{pfflemma}\footnote{There were three typos in Lemma \ref{pfflemma} (all concerning $\mathrm{volume}\left(\Gamma\right)$) that appeared in the first three versions of this paper on the arxiv.  The correct statement and proof appears here.}
Let $\Gamma$ be a compact open subset of $\gln({\Oe}_v)$, and let $F$ be a locally constant Schwartz function
\begin{align*}
F: \Hom_{K_v}(V_v, V_{d, v})\oplus\Hom_{K_v}(V_v, V^d_v)&\rightarrow \bar{\IQ}\\
(X_1, X_2)&\mapsto F(X_1, X_2)
\end{align*}
(with $R$ a subring of $\IC$) whose support in the first variable is $\Gamma$ and such that
\begin{align*}
F(X, { }^tX^{-1}Y) = \chi_1\chi_2^{-1}(\det(X))F(1, Y)
\end{align*}
for all $X$ in $\Gamma$ and $Y$ in $M_n(E_v)$.  Then the Fourier coefficient 
\begin{align*}
c(\beta, 1; f^{PF(-X, Y)}) =: f^{PF(-X, Y)}_{\beta} 
\end{align*} 
of $f^{PF(-X, Y)}$ at $\beta\in M_n(E_v)$  is
\begin{align*}
c(\beta, 1; f^{PF(-X, Y)}) = f^{PF(-X, Y)}_{\beta}(1) = \mathrm{volume}(\Gamma)\cdot F(1, { }^t\beta).
\end{align*}
\end{lem}
\begin{proof}
Throughout the proof, we drop the subscript $v$ from $\Phi$ (to simplify the notation).  For any Schwartz function\footnote{A typo in the exponent of $\left|\det X\right|$ the proof of Lemma \ref{pfflemma} appeared in the first three versions of this paper on the arxiv.} 
\begin{align*}
\Phi: \Hom_{K_v}(V_v, V_{d, v})\oplus\Hom_{K_v}(V_v, V^d_v)\rightarrow \IC,
\end{align*}
\begin{tiny}
\begin{align*}
f_{\beta}^{\Phi}(1) & = \int_{M_n(E_v)}f^\Phi\left(\begin{pmatrix}0&1\\-1&0\end{pmatrix}\begin{pmatrix}1&N\\0&1\end{pmatrix}\right)\e_p(-\tr\beta N)dN\\
& = \int_{GL_n(E_v)}\chi_1^{-1}\chi_2(\det X)\left|\det X\right|^{4s}\int_{M_n(E_v)}\Phi\left(-X, -XN\right)\e^{\left(-\tr(\beta N)\right)}dN d^\times X\\
& = \int_{GL_n(E_v)}\chi_1^{-1}\chi_2(\det X)\left|\det X\right|^{4s}\int_{M_n(E_v)}\Phi\left(-X, N\right)\e^{\left(\tr(\beta X^{-1}N)\right)}dN d^\times X\\
& = \int_{GL_n(E_v)}\chi_1^{-1}\chi_2(\det X)\left|\det X\right|^{4s}P\Phi(-X, -\beta X^{-1}) d^\times X.
\end{align*}
\end{tiny}
So
\begin{align*}
f_{\beta}^{PF(-X, Y)}(1) & = \int_{GL_n(E_v)}\chi_1^{-1}\chi_2(\det X)\left|\det X\right|^{4s}PPF(X, -{ }^tX^{-1}{ }^t\beta) d^\times X\\
& = \int_{GL_n(E_v)}\chi_1^{-1}\chi_2(\det X)\left|\det X\right|^{4s}F(X, { }^tX^{-1}{ }^t\beta) d^\times X\\
& = \int_{\Gamma}F(1, { }^t\beta) d^\times X\\
& = \mathrm{volume}(\Gamma)\cdot F(1, { }^t\beta).
\end{align*}
\end{proof}

\paragraph{A choice of a Schwartz function}
We now choose specific Schwartz functions $F$ meeting the criteria of Lemma \ref{pfflemma}.  These Siegel sections play a key role in our construction of the Eisenstein measure.

Let $(a, b) = ((a(\sigma), b(\sigma))_{\sigma\in\Sigma}$ be the signature of $\langle v_1, v_2\rangle_V$.  For each $v\in\Sigma$, let 
\begin{align*}
a_v & = n_{1, v} + \cdots + n_{t(v), v}\\
b_v & = n_{t(v)+1, v}+ \cdots + n_{r(v), v}
\end{align*}
be partitions of $a(\sigma_v)$ and $b(\sigma_v)$.
 



For each $v\in \Sigma$, let $\mu_{1, v} , \ldots, \mu_{r(v), v}$ be characters of ${\Oe}_v^\times$, and let $\mu_v = (\mu_{1, v}, \ldots, \mu_{r(v), v})$.  Let $\mu = \prod_{v\in\Sigma}\mu_v$.  View each character $\mu_{i,v}$ as a character of $GL_{n_i}\left({\Oe}_v\right)$ via composition with the determinant.  

%

Denote by $\mathfrak{X}$ the subset of $M_n({\Oe}_v)$ consisting of matrices $\begin{pmatrix}A & B\\ C& D\end{pmatrix}$, with $A\in M_{a_v}({\Oe}_v)$, such that the determinant of each of the leading principal $n_{1, v}+\cdots +n_{i, v}$-th minors of $A$ is in ${\Oe}_v^\times$ for $i = 1, \ldots, t(v)$, and the determinant of each of the leading principal $n_{t(v)+1, v}+\cdots +n_{i, v}$-th minors of $D$ is in ${\Oe}_v^\times$ for $i = t(v)+1, \ldots, r(v)$.  

Denote by $A_i$ the determinant of the leading principal $i$-th minor of $A$ and by $D_i$ the leading principal $i$-th minor of $D$.

Let 
\begin{align*}
\nu_{i, v} = \chi_{1, v}^{-1}\chi_{2, v}\mu_{i, v}
\end{align*}
for all $1\leq i \leq r(v)$, and let $\nu_v = \left(\nu_{1, v}, \ldots, \nu_{r(v), v}\right)$.
Define $\phi_{\nu_v}$ to be the function supported on $\mathfrak{X}$ and defined for $X = \begin{pmatrix}A & B\\ C& D\end{pmatrix}\in\mathfrak{X}$ by
\begin{align*}
\phi_{\nu_v}(X) & = \prod_{i = 1}^{t(v)-1} (\nu_{i, v}\cdot \nu_{i+1, v}^{-1})(A_{n_{1, v}+\cdots+ n_{i, v}})\cdot\nu_{t(v), v}(A)\\&\times\prod_{i=t(v)+1}^{r(v)-1}(\nu_{i, v}\cdot \nu_{i+1, v}^{-1})(D_{n_{t(v)+1, v}+\cdots+ n_{i, v}})\cdot \nu_{r(v), v}(D).
\end{align*}

Let
$c\geq \max_{v\in\Sigma, 1\leq i\leq r}(1, \ord_p(\cond(\mu_{i, v}))),$
and let 
$\Gamma = \Gamma(c) = \prod_{v\in\Sigma}\Gamma_v(c)$
be the subgroup of $I = \prod_{v\in\Sigma}I_v$ whose terms below the $n_{i, v}\times n_{i, v}$-blocks along the diagonal are in $\mathfrak{p}_v^c{\Oe}_v$ for all $v$ and such that the upper right $a(v)\times b(v)$ block is also in $\mathfrak{p}_v^c{\Oe}_v$ for all $v$.  For each matrix $m\in\Gamma_v$ with $n_{i, v}\times n_{i,v}$-blocks $m_i$ running down the diagonal, we define
\begin{align*}
\mu_v(m) = \prod_i \mu_{i, v}\left(\det\left(m_i\right)\right).
\end{align*}

Let $\Phi_1$ be the function of $M_{n\times n}(\Oe)$ supported in $\Gamma_v(c)$ (and extended by $0$ to all of $M_{n\times n}({\Oe}_v)$) such that
\begin{align}\label{phimu}
\Phi_{1, v}(x) = \phi_{\mu, v}(x), 
\end{align}
for all $x\in\Gamma_v(c)$.  Let $\Phi_{2, v}$ be the function on $M_{n\times n}\left(E_v\right)$ defined by
\begin{align}\label{Phi2defn}
\Phi_{2, v}(x) = \hat{\phi}_{\nu_v}\left(x\right)
\end{align}
for all $x\in M_{n\times n}\left(E_v\right)$.  In Equation \eqref{Phi2defn}, $\hat{\phi}_{\nu_v}$ denotes the Fourier transform of $\phi_{\nu_v}$, as discussed above.

Define\footnote{Prior versions of this paper had $2X_2$ instead of $X_2$ on the right hand side of Equation \eqref{Phi040812}.  That was a typo.}
\begin{align}\label{Phi040812}
\Phi_{\chi_v, \mu_v} (X_1, X_2) = \mathrm{volume}(\Gamma_v)^{-1}\Phi_{1, v}(-X_1)\cdot \hat{\phi}_{\nu, v}(X_2),
\end{align}
where
\begin{align*}
\nu_v = \chi_{1, v}^{-1}\cdot\chi_{2, v}\cdot \mu_v
\end{align*}
and $\hat\phi_{\nu, v}$ denotes the Fourier transform of the single-variable function $\phi_{\nu, v}$ (so that $\Phi_{\chi_v, \mu_v}$ is a partial Fourier transform in the second variable of a two-variable function, in the sense of Lemma \ref{pfflemma}).
We define
\begin{align}\label{secpequ}
f^{\chi, \mu_v} = f^{\Phi_{\chi, \mu_v}}.
\end{align}
By Lemma \ref{pfflemma}, the $\beta$-th Fourier coefficient of $f^{\chi, \mu_v}$ is
\begin{align}\label{mucoeff}
f^{\chi, \mu_v}_{\beta} = \phi_{\nu, v}({ }^t\beta).
\end{align}
Observe that $f^{\chi, \mu_v}_{\beta} = 0$ whenever $\beta$ is of rank less than $n$.  Also, note that the above definitions and discussion still hold if we replace each character $\mu_v$ by a locally constant function.

We write $\mu = \prod_{v\in\sigma}\mu_v$, $\nu = \prod_{v\in\sigma}\nu_v$, $\phi_{\mu} = \prod_{v\in\Sigma}(\phi_{\mu, v})$, and $f^{\chi, \mu} = \otimes_{v\in\Sigma}f^{\chi, \mu_v}$.
%
%

\begin{rmk}In this remark, we suppress the subscript $v$ where it is clear that we are working at the place $v$.
Let $\delta$ be an $n\times n$ matrix, and let $\gamma$ be an invertible $n\times n$ matrix with entries in $\Oe$.  Note that one can modify the definition of the local sections given in Equation \eqref{sectionininducedrepnatp} so that the integral is again over $\gln({\Oe}_v)$ but $\Phi$ is evaluated at $(X(\gamma C\gamma^{-1}), X (C\delta + D))$ for each $g = \begin{pmatrix}A & B\\ C& D\end{pmatrix}$.  Then we can define $\Phi$ so that we obtain a different section from the one given above, but with similar Fourier coefficients to the ones given above (the main difference being the support of the Fourier coefficients).  There are also other relatively minor modifications one can make to this construction to obtain other local Siegel sections at $p$ with Fourier coefficients of a similar form.  

We note, though, that while these sections all can be used to define Eisenstein series with nice Fourier coefficients at $p$, not all of these coefficients interpolate nicely.
\end{rmk}

\subsubsection{The sections at finite places of $E$ not dividing $p$ or $\infty$}\label{sectionsnotpinfty}
Let $\mathfrak{b}$ be an ideal in $\Oe$ prime to $p$.  For each finite place $v$ prime to $p$, \cite[Section 18]{sh} explains how to define sections $f_v^{\mathfrak{b}} =f_v^{\mathfrak{b}}(\bullet; \chi_v, s)\in\Ind_{P(E_v)}^{G(E_v)}(\chi_v, s)$ with the following property: By \cite[Proposition 19.2]{sh}, whenever the Fourier coefficient $c(\beta, m(1); f_v^{\mathfrak{b}})$ is non-zero,
\begin{tiny}
\begin{align}\label{coeffnotpinfty}
\prod_{v\ndivides p\infty}c(\beta, m(1); f_v^{\mathfrak{b}}) = N_{E/\IQ}(\mathfrak{b}\Oe)^{-n^2}\prod_{i = 0}^{n-1}L^p\left(2s-i, \chi_{E}^{-1}\tau^i\right)^{-1}\prod_{v\ndivides p\infty}P_{\beta, v, \mathfrak{b}}\left(\chi_E(\pi_v)^{-1}\left|\pi_v\right|_v^{2s}\right),
\end{align}
\end{tiny}
where:
\begin{enumerate}
\item{the product is over primes of $E$;}
\item{the Hecke character $\chi_E$ is the restriction of $\chi$ to $E$;}
\item{the function $P_{\beta, v, \mathfrak{b}}$ is a polynomial that is dependent only on $\beta$, $v$, and $\mathfrak{b}$ and has coefficients in $\ZZ$ and constant term $1$;}
\item{the polynomial $P_{\beta, v, \mathfrak{b}}$ is identically $1$ for all but finitely many $v$,}
\item{$\tau$ is the Hecke character of $E$ corresponding to $K/E$,}
\item{$\pi_v$ is a uniformizer of $O_{E, v}$, viewed as an element of $K^\times$ prime to $p$.}
and
\item{
\begin{align*}
L^p(r, \chi_E^{-1}\tau^i) = \prod_{v\ndivides p\infty\mathrm{cond}{\tau}}\left(1-\chi_v(\pi_v)^{-1}\tau^i(\pi_v)\left|\pi_v\right|_v^r\right)^{-1}.
\end{align*}
}
\end{enumerate}
Note that only the factor $\prod_{v\ndivides p\infty}P_{\beta, v}\left(\chi_E(\pi_v)^{-1}\left|\pi_v\right|_v^{2s}\right)$ depends on $\beta$.
\begin{rmk}
The above sections away from $p$, as constructed in \cite{sh}, are built from characteristic functions of lattices (which one can choose to have certain properties corresponding to the choice of ideal $\mathfrak{b}$ and the desired level of the Eisenstein series away from $p$).  These are the same as in \cite{SHL}.
\end{rmk}

\subsubsection{Global Fourier coefficients}

Recall that by Equation \eqref{mavsm1}, the Fourier coefficients $c(\beta, h; f)$ are completely determined by the coefficients $c(\beta, 1_n; f)$.  In Proposition \ref{globalcoeffsprop}, we combine the results of Sections \ref{sectionsinfty}, \ref{sectionatp}, and \ref{sectionsnotpinfty} in order to give the global Fourier coefficients of the Eisenstein series $E_f$.

Let $\chi$ be a unitary Hecke character meeting the conditions of Section \ref{unitaryheckecharacter}, and furthermore, suppose the infinity type of $\chi$ is 
\begin{align}
\prod_{\sigma\in\Sigma}\sigma^{-k-2\nu(\sigma)}\left(\sigma\bar{\sigma}\right)^{\frac{k}{2}+\nu(\sigma)}\label{integralk}
\end{align}
(i.e. $k(\sigma) = k\in\ZZ$ for all $\sigma\in\Sigma$).  Let $c(n, K)$ be the constant dependent only upon $n$ and $K$ defined in Equation \eqref{CnK}.

\begin{prop}\label{globalcoeffsprop}
Let $k\geq n$, and let
\begin{align}\label{aholosec}
f_{k, \nu, \chi, F} := f_{k, \nu, \chi, \mathfrak{b}, F} := \otimes_{v\in\Sigma}f_{F,v}\otimes f^{k, \nu}_{\infty}\left(\bullet; i1_n, \chi, \frac{k}{2}\right)\otimes f^{\mathfrak{b}}\in \Ind_{P\left(\adeles_E\right)}^{G\left(\adeles_E\right)}\left(\chi\cdot\left|\cdot\right|_K^{-\frac{k}{2}} \right),
\end{align}
with $\chi$ as in Equation \eqref{integralk}, $\otimes_{v\divides\Sigma}f_{\chi, \mu}$ the section at $p$ defined in Equation \eqref{secpequ}, $f^{k, \nu}_{\infty}$ the section at $\infty$ defined in Section \ref{sectionsinfty}, and $f^{\mathfrak{b}}$ the section away from $p$ and $\infty$ defined in Section \ref{sectionsnotpinfty}.
 
Then all the nonzero Fourier coefficients $c(\beta, 1_n; f_{k, \nu, \chi, F})$ are given by
\begin{tiny}
\begin{align}\label{FourierCoeffFormula}
 D(n, K, \mathfrak{b}, p, k)\prod_{v\ndivides p\infty}P_{\beta, v, \mathfrak{b}}\left(\chi_E(\pi_v)^{-1}\left|\pi_v\right|_v^{k}\right)\phi_{\chi_1^{-1}\chi_2\mu}\left({ }^t\beta\right)\prod_{v\in\Sigma}\sigma_v(\det\beta)^{k-n}\e\left(i\tr_{E/\IQ}(\beta)\right).
 \end{align}
\end{tiny}
where
\begin{tiny} 
\begin{align*}
D&(n, K, \mathfrak{b}, p, k)\\  
&= C(n, K)N(\mathfrak{b}{\Oe})^{-n^2}\prod_{i = 0}^{n-1}\left(2^{(1-n)n}i^{-nk}(2\pi)^{nk}\left(\pi^{n(n-1)/2}\prod_{t=0}^{n-1}\Gamma(k-t)\right)^{-1}\right)^{[E:\IQ]}\prod_{i = 0}^{n-1}L^p\left(k-i, \chi_{E}^{-1}\tau^i\right)^{-1}.
\end{align*}
\end{tiny}

\end{prop}

\begin{proof}
This follows directly from Section \ref{consequencesAB}, Lemma \ref{pfflemma}, and Equations \eqref{coeffnotpinfty} and \eqref{coeffinfty}.  
\end{proof}

Let 
\begin{align*}
\alpha(\beta) := \alpha(\beta, \mathfrak{b}, \chi):=
\prod_{v\ndivides p\infty}P_{\beta, v, \mathfrak{b}}\left(\chi_E(\pi_v)^{-1}\left|\pi_v\right|_v^{k}\right).
\end{align*}
Note that $\alpha(\beta, \chi)$ is a (finite) $\ZZ$-linear combination of terms of the form $\prod_{v\ndivides p\infty}\chi_v(a)^{-1}\left|a\right|_v^k$, with $a$ a $p$-integral element of the integer ring of the totally real field $E$.  For any $p$-integral element $a$ of the totally real field $E$, we have
\begin{align}
\prod_{v\ndivides p\infty}\chi_v(a)^{-1}\left|a\right|_v^k &= \chi_1\chi_2^{-1}(a)\prod_{v\in\Sigma}\sigma_v(a)^{-k-2\nu(\sigma)}\left(\sigma_v(a)\bar{\sigma}_v(a)\right)^{\frac{k}{2}+\nu(\sigma)}\left(\prod_{v\in\Sigma}\sigma_v(a)\bar\sigma_v(a)\right)^{-k/2}\nonumber\\
&= \chi_1\chi_2^{-1}(a)\prod_{v\in\Sigma}\sigma_v(a)^{-k}.\label{ridofalpha}
\end{align}

\subsubsection{$q$-expansions}\label{qexpansions}

Let $f = \otimes_vf_v\in \otimes_vI(\chi_v, s)$ be a section whose component at $\infty$ is $f_{\infty}^{k, \nu}$, defined as in Section \ref{sectionsinfty}.
For $z = x+iy\in \hn$, $h\in GL_n(\adeles_{K, f})$, and $\alpha\in \prod_{v\divides\infty} G_v$ such that $\alpha \gimel = z$, we put
\begin{align*}
E_{k, \nu,}(z; h, \chi, \mu, s) := j_{\alpha}^{(k(\sigma)_{\sigma}, \nu(\sigma)_{\sigma})}(\gimel) E_f(m(h) \alpha).
\end{align*}
Then, as explained in \cite[Lemma 18.7(2)]{sh}, one can give the Fourier expansion for the function $E_f(z; h, k, \nu, s)$ of the variable $z$ in terms of the Fourier coefficients given above.  (Note that Equation \eqref{gimelvsi} gives the correction factor for using $\gimel\in\hn$ instead of $i$.)

In particular, when the infinity type of $\chi$ is as in \eqref{integralk} and $f$ is the Siegel section defined in Equation \eqref{aholosec}, $E_{k, \nu}\left(z; h, \chi, \mu, \frac{k}{2}\right)$ is a holomorphic function of $z$ such that
\begin{small}
\begin{align*}
&D(n, K, \mathfrak{b}, p, k) ^{-1}E_{k, \nu}\left(z; 1, \chi, \mu, \frac{k}{2}\right)  =\\
& \sum_{\beta} \prod_{v\ndivides p\infty}P_{\beta, v, \mathfrak{b}}\left(\chi_E(\pi_v)^{-1}\left|\pi_v\right|_v^{k}\right)\phi_{\chi_1^{-1}\chi_2\mu}\left({ }^t\beta\right)\prod_{v\in\Sigma}\sigma_v(\det\beta)^{k-n}\e_{\infty}{(\tr(\beta z))}.
\end{align*}
\end{small}
(The above sum is over $\beta\in\hern$ such that the Fourier coefficient at $\beta$ is nonzero.)  Observe that for this choice of $f$, the Fourier coefficients of $E_f(z; h)$ are independent of $z$ and algebraic.  Furthermore, if we choose $F$ so that $F(1, \beta)$ is $p$-integral for every $\beta$, then the Fourier coefficients of $E_f(z; h)$ are also $p$-integral.  Let $\mathcal{K}$ be the compact subgroup of $G(\adeles)$ so that we may view $E_{k, \nu, F}$ as an automorphic form on the complex points of $_{\mathcal{K}}Sh(W)$ (i.e. $\mathcal{K}$ corresponds to the above choice of local Siegel sections).


Denote by $G_{k, \nu, \chi, \mu}$ the automorphic form on the $_{\mathcal{K}}Sh(W)$ whose algebraically defined $q$-expansion (i.e. its value at the Mumford object, generalizing the Tate curve; see \cite{la, lanalgan, ha86} for more on algebraic $q$-expansions) at a cusp $L$ is 
%

\begin{align}\label{Gknuchimu}
G_{k, \nu, \chi, \mu}(q) = \sum_{\beta\in L}\left(\alpha(\beta, \chi)\phi_{\chi_1^{-1}\chi_2\mu}\left({ }^t\beta\right)\prod_{v\in\Sigma}\sigma_v(\det\beta)^{k-n}\right)q^{\beta},
\end{align}
where
\begin{align*}
\alpha(\beta, \chi)=
\prod_{v\ndivides p\infty}P_{\beta, v, \mathfrak{b}}\left(\chi_E(\pi_v)^{-1}\left|\pi_v\right|_v^{k}\right),
\end{align*}
$k>2n-1$, and $\chi$ is a Hecke character meeting the conditions of Section \ref{eisensteinseriesprelims}.
Note that the $q$-expansion coefficients of $G_{k, \nu, \chi, \mu}(q)$ are the same as those in the Fourier expansion of the holomorphic function $D(n, K, \mathfrak{b}, p, k)^{-1}E_{k, \nu}(z; 1, \chi, \frac{k}{2})$ given above, and the value of $G_{k, \nu, \chi, \mu}$ at a $\IC$-valued point of $_{\mathcal{K}}Sh(W)$ is the same as the value of $E_{k, \nu}$ at the corresponding point $z$.

An argument similar to \cite[Theorem (3.4.1)]{kaCM} shows that $G_{k, \nu, \chi, \mu}(q)$ is, in fact, the $q$-expansion of a $p$-adic automorphic form (simply the image under the canonical map from the space of automorphic forms over $\OCp$ to the space of $p$-adic automorphic forms), which we also denote by $G_{k, \nu, \chi, \mu}$.

\begin{rmk}{Relationship with Katz's Eisenstein series for CM-fields}


Plugging $n=1$ into Equation \eqref{FourierCoeffFormula} and comparing with \cite[Equations (3.2.6)-(3.2.8)]{kaCM} shows that for $n=1$, our Eisenstein series are similar to those in \cite{kaCM}, and for $n>1$, our Eisenstein series (and the Fourier coefficients in Equation \eqref{FourierCoeffFormula}) are a natural generalization of those in \cite{kaCM}.

In particular, since the terms in Equation \eqref{FourierCoeffFormula} look more complicated that in \cite{kaCM}, we note that due to cancellation when $n=1$, we obtain
\begin{align*}
c(1, K) &= \left|D_E\right|^{-1/2},\\
D(1, K, \mathfrak{b}, p, k) &= \left|D_E\right|^{-1/2}N(\mathfrak{b})\left(\left(-2\pi i\right)^k\Gamma(k)^{-1}\right)^{[E:\IQ]}L^p(k, \chi_E^{-1})^{-1}.
\end{align*}
\end{rmk}

\subsubsection{Certain $p$-adic Eisenstein series}\label{certaineseries}
Via the isomorphism \eqref{okoeid}, we identify $\OK\otimes\ZZ_p$ with $\left(\Oe\otimes\ZZ_p\right)\times\left(\Oe\otimes\ZZ_p\right)$.  For each $v\in\Sigma$, let $r_v = r(v)$ be a positive integer, and let $r = \left(r_v\right)_v$.  Let
\begin{align}\label{Tofr}
T(r) = \prod_{v\in\Sigma}{\underbrace{{\Oe}_v^\times\times\cdots \times{\Oe}_v^\times}_{r_v \mbox{ copies}}}.
\end{align}

Let $\mu$ be a locally constant function on $T(r)$, extended by $0$ to all of $\prod_{v\in\Sigma}{\underbrace{{\Oe}_v\times\cdots \times{\Oe}_v}_{r_v \mbox{ copies}}}$.  Fix integers $k$ and $\nu(\sigma)$, and let $\nu = (\nu(\sigma))_{\sigma}$.  Let $F$ be a locally constant function on $\left(\OK\otimes\ZZ_p\right)^{\times}$, extended by $0$ to all of $\OK\otimes\ZZ_p$, such that
\begin{align}\label{Fcond}
F(e x) = \prod_{\sigma\in\Sigma}\sigma(e)^{k+2\nu(\sigma)}\left(\sigma(e)\overline{\sigma}(e)\right)^{-\left(\frac{k}{2}+\nu(\sigma)\right)}F(x)
\end{align}
for all $e\in\OK^{\times}$ and $x\in\left(\OK\otimes\ZZ_p\right)^{\times}$.  Then there is an integer $d$, elements $c_1, \ldots, c_d\in \OK$, and Hecke characters $\xi_1, \ldots, \xi_d$ of $K$ of conductor dividing $p^{\infty}$ and infinity type $\prod_{\sigma\in\Sigma}\sigma^{k+2\nu(\sigma)}\left(\sigma\overline{\sigma}\right)^{-\left(\frac{k}{2}+\nu(\sigma)\right)}$ such that
\begin{align*}
F(a) &= c_1F_1(a) + \cdots + c_dF_d(a),\\
F_i &= \prod_{v\divides p}\xi_{i, v}
\end{align*}
for all $a\in \left(\OK\otimes\ZZ_p\right)^{\times}$.  Let
\begin{align*}
\psi_{k, \nu} = \prod_{\sigma\in\Sigma}\sigma^{-k}\left(\frac{\bar{\sigma}}{{\sigma}}\right)^{\nu(\sigma)}.
\end{align*}
\begin{lem}\label{lclem}
Given such a locally constant function $F$, there is an automorphic form $E_{F, \mu}$ defined over $\OK$ of weight $k, \nu$, whose $q$-expansion is given by
\begin{tiny}
\begin{align*}
E_{F, \mu}(q) &= \sum_{i=1}^d E_{F_i, \mu}(q)\\
E_{F_i, \mu}(q) &= \sum_{\beta\in L}c_i\left(\left(\prod_{v\ndivides p\infty}P_{\beta, v, \mathfrak{b}}\left(F_i(\pi_v)\psi_{k, \nu}\left(\pi_v\right)\right)\right)\cdot F_i(\det\beta^{-1})\psi_{k, \nu}(\det\beta^{-1})\prod_{v\in\Sigma}\sigma_v\left(\det\beta\right)^{-n}\phi_\mu(\beta)\right)q^{\beta}.
\end{align*}
\end{tiny}
\end{lem}
\begin{proof}
We take $E_{F_i, \mu}$ to be the automorphic form $G_{k, \nu, \xi_i, \mu}$ in Equation \eqref{Gknuchimu}.  The lemma then follows immediately.
\end{proof}
\begin{rmk}
Similarly to the discussion in \cite[Section 3]{kaCM}, we may generalize the statement of Lemma \ref{lclem} to the case of locally constant functions $F$ and $\mu$ having image in any $\OK$-algebra $R$ (in which case $E_{F, \mu}$ is an automorphic form over $R$).
\end{rmk}

Suppose now that $F$ is a continuous function satisfying Equation \eqref{Fcond} supported on $\left(\OK\otimes\ZZ_p\right)^\times$ and that $\mu$ is a continuous function on $T(r)$ extended by $0$ to all of $\prod_{v\in\Sigma}{\underbrace{{\Oe}_v\times\cdots \times{\Oe}_v}_{r_v \mbox{ copies}}}$.  Then there is a $p$-adic automorphic form $E_{F, \mu}$ obtained via the $q$-expansion principle and an argument similar to the proof of \cite[Theorem (3.4.1)]{kaCM}.  The $q$-expansion coefficients of $E_{F, \mu}$ are essentially $p$-adic limits of $q$-expansion coefficients of automorphic forms $E_{F_i, \mu_i}$ with $F_i$, $\mu_i$ locally constant and $F_i$ satisfying Equation \eqref{Fcond}.  When $\mu$ and $F$ are locally constant, this automorphic form is defined over $\OK$ and is just the image of the automorphic form from Lemma \ref{lclem} under the canonical map from the space $\mathcal{A}(k, \nu, \mu, \OK)$ of weight $k, \nu$ automorphic forms on $_{\mathcal{K}}Sh(W)$ into the space $\mathcal{V}$ of $p$-adic automorphic forms (viewed as sections on the Igusa tower).

Thus, it follows that for {\it any} continuous functions $F$ supported on $\left(\OK\otimes\ZZ_p\right)^\times$ satisfying Equation \eqref{Fcond} (for {\it some} integers $k$ and $\nu$) and $\mu$ supported on $\prod_{v\in\Sigma}{\underbrace{{\Oe}_v^\times\times\cdots \times{\Oe}_v^\times}_{r_v \mbox{ copies}}}$ with values in an $\OK$-algebra $R$, we obtain a $p$-adic automorphic form $E_{F, \mu}$ defined over $R$ whose $q$-expansion coefficients are limits of $q$-expansion coefficients of some automorphic forms $E_{F_i, \mu_i}$ for some locally constant $F_i$'s and $\mu_i$'s.

\begin{rmk}
When $F$ and $\mu$ are locally constant $\bar{\IQ}$-valued functions, $E_{F, \mu}$ is just a sum of the analytically defined Eisenstein series discussed earlier.
\end{rmk}

\subsubsection{Hecke characters of type $A_0$, viewed $p$-adically}\label{viewedpadically}

Let
\begin{align*}
\chi: K^\times\backslash\idelesK\rightarrow \IC^{\times}
\end{align*}
be a Hecke character of type $A_0$, so on  $K$, 
\begin{align*}
\chi_{\infty} (a)= \iota_{\infty}\circ\left(\prod_{\sigma\in\Sigma}\lp\frac{1}{\sigma}\rp^k\lp\frac{\bar{\sigma}}{\sigma}\rp^{d(\sigma)}\right)
\end{align*}
with $k$ and $d(\sigma)$ in $\ZZ$ for all $\sigma$.  Note that on $K$,
\begin{align*}
\iota_p\circ\left(\prod_{\sigma\in\Sigma}\lp\frac{1}{\sigma}\rp^k\lp\frac{\bar{\sigma}}{\sigma}\rp^{d(\sigma)}\right)
\end{align*}
defines a $p$-adic character that extends continuously to a ($p$-adically) continuous $\bar{\IQ}_p$-valued character on $K\otimes\IQ_p$.  So we can define a $p$-adic Hecke character on the finite id\`eles
\begin{align}\label{tildechi}
\tilde{\chi}=\prod_{v\ndivides \infty}\tilde{\chi}_v:K^\times\backslash{\idelesK}^{,\infty}\rightarrow\bar{\IQ}_p^\times
\end{align}
by
\begin{align*}
\tilde{\chi}_v & = \chi_v, \mbox{ if $v\ndivides p$}\\
\tilde{\chi}_p & = \chi_p\cdot \left(\iota_p\circ\left(\prod_{\sigma\in\Sigma}\lp\frac{1}{\sigma}\rp^k\lp\frac{\bar{\sigma}}{\sigma}\rp^{d(\sigma)}\right)\right).
\end{align*}
Note that the restriction of $\tilde{\chi}$ to $\adeles_K^{\times, p, \infty}\times \lp\OK\times\ZZ_p\rp^\times$ gives a $p$-adic character
\begin{align*}
\adeles_K^{\times, p, \infty}\times\lp\OK\times\ZZ_p\rp^\times\rightarrow \OCp^\times.
\end{align*}

When we make the shift
\begin{align*}
\chi\rightarrow\tilde{\chi},
\end{align*}
we shall refer to ``shifting to $p$."

\begin{example}\label{Fexl}  Let $\nu = (\nu(\sigma))_{\sigma\in\Sigma}\in\ZZ^{\Sigma}$.  Let $\chi$ be a Hecke character of type $A_0$ and whose conductor divides $p^{\infty}$ and whose infinity-type is
\begin{align*}
\psi_{k, \nu} = \prod_{\sigma\in\Sigma}\sigma^{-k}\left(\frac{\bar{\sigma}}{{\sigma}}\right)^{\nu(\sigma)}.
\end{align*}
Let
\begin{align*}
F = \prod_{v\divides p}\tilde{\chi}_v,
\end{align*}
where $\tilde{\chi}$ is defined as in Equation \eqref{tildechi}.
Then $F(e) = 1$ for all $e\in \OK^{\times}$, and
\begin{align}\label{EFmu}
E_{F, \mu}(q) & = \sum_{\beta}\left(\left(\prod_{v\ndivides p\infty}P_{\beta, v, \mathfrak{b}}\left(F(\pi_v\right))\right)\cdot F(\det\beta^{-1})N_{E/\IQ}(\det\beta)^{-n}\phi_\mu(\beta)\right)q^{\beta}.
\end{align}
If $\chi' = \chi\left|\cdot\right|_{K/\IQ}^{\frac{k}{2}}$, then
\begin{align*}
E_{F, \mu} = G_{k, \nu, \chi', \mu}.
\end{align*}
\end{example}

\subsection{Some comments about a functional equation and the connection with Katz's setting}\label{fcnlequ}
I would like to thank the anonymous referee for suggesting that I include some comments about a functional equation. We provide some brief comments on this  topic here.
%

As above, let $\chi$ be a unitary Hecke character whose infinity type is as in Equation \eqref{integralk} and whose conductor divides $p^\infty$.  

For simplicity of notation and to remain as close as possible to Shimura's setup, we omitted similitude factors above.  To view the functional equation from \cite{kaCM} in terms of our Eisenstein series, however, it is helpful to introduce similitude factors.  Let $\nu: GU\left(\eta, \adeles_E\right)\rightarrow \adeles_E^\times$ denote the similitude character.  We extend the Siegel sections from above to Siegel sections $f\in\Ind_{GP}^{GU}\left(\chi, s\right)$ so that $f\left(\diag\left({ }^t\bar{h}^{-1}, \lambda h\right)\right) = \chi\left(\det\left(\lambda h\right)\right)|\det(\lambda h)|^{-s}$ for all $\lambda\in \adeles_E^\times$.  (Note that $\lambda = \nu\left(\diag\left({ }^t\bar{h}^{-1}, \lambda h\right)\right)$.)

%

Lemma \ref{mavsm1} then becomes 
\begin{lem}
For each $h\in \gln(\adeles_K)$, $\lambda\in\adeles_E^\times$, and $\beta\in\hern(K)$,
\begin{align}\label{fcoeffsat1}
c\left(\beta, \begin{pmatrix}{ }^t\bar{h}^{-1}& 0\\ 0& \lambda h\end{pmatrix}; f\right) = \left|\lambda\right|^{n^2}\chi({ }^t\bar{h}^{-1})\left|\det\left({ }^t\bar{h}^{-1}\cdot h^{-1}\right)\right|_E^{n-s} c(\lambda^{-1}{ }^t\bar{h}^{-1}\beta h^{-1}, 1_n; f).
\end{align}
\end{lem}
\begin{proof}
Let $\eta = \begin{pmatrix}0 & -1_n\\ 1_n & 0\end{pmatrix}.$
Let $m(h, \lambda)$ denote the matrix $\begin{pmatrix}{ }^t\bar{h}^{-1} & 0\\ 0 & \lambda h\end{pmatrix}$.  Observe that for any $n\times n$ matrix $m$,
\begin{align*}
\eta \cdot m(h, \lambda)\cdot \eta^{-1} & = m(\lambda^{-1}{ }^t\bar{h}^{-1}, \lambda)\\
m(h, \lambda)^{-1}\cdot \begin{pmatrix}1&m\\ 0 &1\end{pmatrix}\cdot m(h, \lambda) &= \begin{pmatrix}1 & \lambda { }^t\bar{h} m h\\ 0 & 1\end{pmatrix}.
\end{align*}
Therefore,
\begin{align*}
\eta \cdot \begin{pmatrix}1 & m \\ 0 & 1\end{pmatrix}\cdot m(h, \lambda) & = \left(\eta \cdot m(h, \lambda)\cdot \eta^{-1}\right)\eta \left(m(h, \lambda)^{-1}\begin{pmatrix}1 & m\\ 0 & 1\end{pmatrix} m(h, \lambda)\right)\\
 & = m(\lambda^{-1}{ }^t\bar{h}^{-1}, \lambda) \eta \begin{pmatrix}1 & \lambda{ }^t\bar{h} m§ h\\ 0 & 1\end{pmatrix}.
\end{align*}
So for any place $v$ of $E$ and section $f_v\in \Ind_{P(E_v)}^{G(E_v)}(\chi, s)$,
\begin{align}\label{leviunimp2}
f_v\left(\eta\begin{pmatrix}1 & m \\ 0 & 1\end{pmatrix}m(h_v, \lambda)\right) = \chi_v\left({ }^t\bar{h_v}^{-1}\right)\left|{ }^t\bar{h_v}^{-1}\right|^{-2s}_vf_v\left(\eta \begin{pmatrix}1 & \lambda { }^t\bar{h}_v m h_v\\ 0 & 1\end{pmatrix}\right). \end{align}
The lemma now follows from Equation \eqref{leviunimp2} and the fact that the Haar measure $d_v$ satisfies $d_v(\lambda h_vx{ }^t\bar{h_v}) = \left|\det\left(\lambda{ }^t\bar{h_v}\cdot h_v\right)\right|_v^{n}d_v(x)$ for each place $v$ of $E$. 
\end{proof}

So
\begin{align*}
c\left(\beta, \begin{pmatrix}\lambda^{-1}{ }^t\bar{h}^{-1}& 0\\ 0& h\end{pmatrix}; f\right) 
& = \chi(\lambda^{-n})\left|\lambda^{-n}\right|^{-s}c\left(\beta, \begin{pmatrix}{ }^t\bar{h}^{-1}& 0\\ 0& \lambda h\end{pmatrix}; f\right)
\end{align*}

Now, the points in $M(\adeles_E)$ parametrize cusps. We have 
\begin{align}
Tate_{{ }^t\bar{h}^{-1}, \lambda h}(q)\dual = Tate_{\lambda^{-1}{ }^t\bar{h}^{-1}, h}(q).\label{tatefe}
\end{align}
So we obtain the following functional equation
\begin{align*}
E(Tate_{{ }^t\bar{h}^{-1}, \lambda h}(q)\dual) &= E(Tate_{\lambda^{-1}{ }^t\bar{h}^{-1}, h}(q))\\
& = \chi^{-1}(\lambda^n)\left|\lambda^{-n}\right|^{-s}E(Tate_{{ }^t\bar{h}^{-1}, \lambda h}(q))
\end{align*}

For $n=1$, we rephrase the above discussion in the language of \cite{kaCM} to relate the functional equation in 
Equation \eqref{functionalequation} to the functional equation in \cite[Section 3]{kaCM}.  Let $\mathfrak{c}$ be a fractional ideal in $E$, and let $\mathfrak{a}$ and $\mathfrak{b}$ be fractional ideals of $K$ such that $\mathfrak{c} = \mathfrak{a}\bar{\mathfrak{b}}^{-1}$.  Let $\lambda$ be an idele class in $\adeles_E^\times$ corresponding to the fractional ideal $\mathfrak{c}^{-1}$.  Let $h$ be an idele class in $\adeles_K^\times$ corresponding to the fractional ideal $\bar{\mathfrak{a}}$.  Then $\lambda h$ is an idele class corresponding to the fractional ideal $\mathfrak{b}$.  Then the abelian variety that \cite{kaCM} denotes by $Tate_{\mathfrak{a}, \mathfrak{b}}(q)$ is the abelian variety that we denote above by $Tate_{{ }^t\bar{h}^{-1}, \lambda h}(q)$.  Now, $Tate_{\mathfrak{a}, \mathfrak{b}}(q)$ is a $\mathfrak{c}$-polarized abelian variety whose dual is the $\mathfrak{c}^{-1}$-polarized abelian variety $Tate_{\bar{\mathfrak{b}}, \bar{\mathfrak{a}}}(q)$ (which in our notation above is $Tate_{\lambda^{-1}{ }^t\bar{h}^{-1}, h}(q)$).  The $q$-expansion sum in then over 
\begin{align*}
\left\{\beta\in \mathfrak{ab}| \beta \mbox{ is totally real}\right\}.  
\end{align*}
For any $\lambda$-polarized Hilbert-Blumenthal abelian variety $\underline{X}$, we have the following functional equation
\begin{align}\label{functionalequation}
E(\underline{X}\dual) = \chi(\lambda^{-n})\left|\lambda^{-n}\right|^{-s}E(\underline{X}).
\end{align}

Once we have expressed our functions $F$ in a form similar to the form of the functions denoted by $F$ in \cite{kaCM}, we shall obtain a functional equation similar to the one in \cite{kaCM} for $n=1$.
As part of our discussion, we express the two-variable locally constant functions that Katz denotes by $F(x, y)$ throughout \cite{kaCM} (first appearing in \cite[Theorem 3.2.3]{kaCM}) explicitly in terms of the two-variable locally constant functions we also denote by $F(x, y)$ in (first appearing in Lemma \ref{pfflemma}).  


Let $\mathfrak{a}$ and $\mathfrak{b}$ be fractional ideals in $K$. Let $\beta = ab\in \mathfrak{ab}$ with $b$ a $p$-adic unit in $\mathfrak{b}$ and $a$ a $p$-adic unit in $\mathfrak{a}$, and furthermore, suppose that $b$ and $\beta$ are totally real and that $\sigma(\beta)\geq 0$ for all $\sigma\in\Sigma$.  Then for any function $F$ satisfying the hypotheses of Lemma \ref{pfflemma} with $n=1$ and $\Gamma = \prod_{v\in\Sigma}{\Oe}_v^\times$,
\begin{align*}
F(1, \beta) &= \chi_1\chi_2^{-1}(b^{-1})F(b, a).
\end{align*}
Note that if $e$ is an element of $\Oe^\times$, 
\begin{align*}
\chi_1\chi_2^{-1}(e^{-1})&=\prod_{\sigma\in\Sigma}\sigma(e^{-1})^{-k}\\
&=\mathbb{N}_{E/\IQ}(e)^k.
\end{align*}
So for any element $e\in\Oe^\times$, 
\begin{align}
F(e^{-1}x, ey)& = \chi_1\chi_2^{-1}(e^{-1})F(x, y)\nonumber\\
&= \mathbb{N}_{E/\IQ}(e)^k F(x, y)\label{katznormcond}
\end{align}
for all $x$ and $y$.  On the other hand, any locally constant function $F$ satisfying Equation \eqref{katznormcond} is a linear combination of functions $F$ of the type we have been considering all along (i.e. satisfying the conditions of Lemma \ref{pfflemma}).  In other words, Katz's function $F$ can be obtained from ours.

Note that
\begin{align}
\Sgn\left(\mathbb{N}_{E/\IQ}(a)\right)\mathbb{N}_{E/\IQ}\left(a\right)^{k-1}F(b, a) &=  \Sgn\left(\prod_{\sigma\in\Sigma}\sigma\left(b^{-1}\beta\right)\right)\prod_{\sigma\in\Sigma}\sigma\left(b^{-1}\beta\right)^{k-1}\prod_{v\ndivides p}\chi_v(b^{-1}) F(1, \beta)\label{begfe}\\
&= \prod_{\sigma\in\Sigma}\sigma(\beta)^{k-1}\prod_{v \ndivides p\infty}\left(\chi_v(b^{-1})\left|b\right|_v^{k-1}\right)F(1, \beta)\nonumber\\
&=\prod_{\sigma\in\Sigma}\sigma(\beta)^{k-1}\prod_{v \ndivides p\infty}\left(\left|\pi_v\right|_v^{-1}\chi_v(\pi_v)^{-1}\left|\pi_v\right|_v^{k}\right)^{r_v}F(1, \beta),\nonumber
\end{align}
where $r_v$ is the integer such that $\left| b\right|_v=\left|\pi_v\right|_v^{r_v}$ (and $\pi_v$ denotes a uniformizer at the place $v$ of $E$).  The product of the polynomials $P_{\beta, v, \mathfrak{b}}$ expresses the Fourier coefficient at $\beta$ as a sum over totally real $b\in\mathfrak{b}$ such that $\beta = ab$ for some $a\in\mathfrak{a}$ (modulo $\Oe^\times$).  We shall not make these polynomials precise here.  However, we note that in the case where $n=1$ and $F$ is identically $1$ (on the support of $F$, with the support of $F$ as above), \cite[Section 18.6]{shar} explains precisely how to derive a formula of the form in \cite{kaCM} from the one in this paper.  The sum is then over totally real $a\in \mathfrak{a}$ and $b\in\mathfrak{b}$ such that $ab = \beta$.  Note that multiplying the left hand side of Equation \eqref{begfe} by $\chi^{-1}(\lambda)|\lambda|^{k-1}$ swaps the roles of $a$ and $b$, as in \cite{kaCM}.  Note that for this conversation, we have omitted extra factors that appear in the Fourier expansion, such as norms of ideals (and id\`eles).  From the perspective of functional equations and writing down Fourier coefficients that look like the ones in \cite{kaCM}, the series in \cite[Equation (18.31.3)]{sh} is perhaps more useful to use in place of the polynomials, because of its more explicit from in Shimura.  From the point of view of $p$-adic interpolation (our main goal in this paper), though, it was more natural to work with the polynomials.

When $n>1$, it is not apparent that there is a similar form of the functional equation (i.e. swapping the entries of $F$).  A close study of the polynomials $P_{v, \mathfrak{b}, \beta}$ (which have not been computed anywhere for $n>1$) or carefully working through \cite[Equation (18.13.3)]{sh} might provide insight, but it seems that for $n>1$, the cleanest form of a functional equation is the one provided in Equation \eqref{tatefe}.

\section{Differential Operators}\label{diffopsnew}
Throughout this section, let $k\in\ZZ$, and let $d = (d(\sigma))_{\sigma\in\Sigma}, \nu = (\nu(\sigma))_{\sigma\in\Sigma}$ with $d(\sigma), \nu(\sigma)\in\ZZ$, $d(\sigma)>0$, for all $\sigma$.  We write $k-d, \nu-d$ to denote $(k-d(\sigma))_{\sigma\in\Sigma}, (\nu(\sigma)-d(\sigma))_{\sigma\in\Sigma}$.  

We briefly review the most important details of certain differential operators that will allow us to $p$-adically interpolate special values of the Eisenstein series $E_{k, \nu}$ (normalized by a period), including certain cases where the Eisenstein series is non-holomorphic (for example, $E_{k, \nu}(z; h, \chi, \mu, s)$ with $s\neq \frac{k}{2}$, as well as the case where $k$ is no longer an integer).  

Since these differential operators are discussed in thorough detail in \cite{kaCM, shclassical, sh, shar, EDiffOps}, we use this section of the paper to summarize only the basic properties necessary for this paper.

\paragraph{$\ci$ differential operators}
There is a differential operator $D_{k, \nu}^d$ that acts on $\ci$-automorphic forms on $\hn$ and satisfies:
\begin{align*}
D_{k, \nu}^{d} \left( j_{\alpha}^{k, \nu}(z)^{-1} \right) = \Psi_{n, d, k}\cdot \det(z-{ }^t\bar{z})^{-d}\cdot j_{\alpha}^{k+2d, \nu-d}(z)^{-1}\left|j_{\alpha}(z)\right|^{2d},
\end{align*}
where
\begin{align*}
\Psi_{n, d, k} = \prod_{\sigma\in\Sigma}\prod_{i = 1}^n\prod_{j=1}^{d(\sigma)}\left(i-j-k\right).
\end{align*}
In particular, 
\begin{align}\label{weightup}
D_{k, \nu}^d \left(E_{k, \nu}\left(z; h, \chi, \mu, \frac{k}{2}\right)\right) &= (\det(z-{ }^t\bar{z})^{-d}\cdot \Psi_{n, d, k})E_{k+2d, \nu-d}\left(z; h, \chi, \mu, \frac{k}{2}\right)
\end{align}
(where $E_{k, \nu}$ is defined as in Section \ref{qexpansions}).

The operator $D_{k, \nu}^d$ is a special case of the $\ci$-differential operators (often called ``Shimura-Maass" operators) studied extensively by Shimura \cite{shclassical, sh, shar} (and also studied algebro-geometrically in \cite{EDiffOps}).  It is the operator that would be denoted $D_{\det^{k+\nu}\otimes\det^{-\nu}}^{\det^d}$ in the conventions of \cite[Equation (12.17)]{shar} (generalized from $SU(\eta)$ to $U(\eta)$, since \cite{shar} only deals with special unitary groups) and \cite[Section 6]{EDiffOps}.\footnote{Note that there should be a $2\pi i$ on the right-hand side of \cite[Equation (3.29)]{EDiffOps}.  The reason for this is explained at the beginning of \cite[Section 4.4]{hasv}.  Correspondingly, the $\ci$-differential operator in \cite{EDiffOps} is actually a power of $2\pi i$ times Shimura's differential operators discussed in \cite{shar, sh}.}
  It is also the operator that would be denoted $D_k(\det^d)$ in the conventions of \cite[Lemma 23.5]{sh}, and it is the generalization to $U(\eta_n)$ of the operator $\Delta_k^d$ in \cite[Equation (17.20)]{shar}.
  
\subparagraph{Differential Operators and Pullbacks}
As demonstrated in Equation \eqref{weightup}, these differential operators change the weight of the automorphic forms to which they are applied.  Correspondingly, they also change the weight of the pullback of such an automorphic form to $U(V)\times U(V)$.  A precise description of the affect of differential operators on the automorphy factors of pullbacks is given in \cite[Sections 23.6-23.9]{sh}.

\paragraph{$p$-adic differential operators and their action on $q$-expansions}
By Theorem IX.3 in \cite{EDiffOps} applied in the special case of scalar-valued automorphic forms, there is a $p$-adic differential operator $\theta^d$ that acts on $p$-adic automorphic forms (viewed as functions on ordinary abelian varieties with PEL structure corresponding to the choice of Shimura data and a particular moduli problem, i.e. as sections over the corresponding Igusa tower, as discussed in \cite{SHL} and \cite{EDiffOps}) 
whose action on the $q$-expansion of a $p$-adic automorphic form $f$ with $q$-expansion
$f(q) = \sum_{\beta} a(\beta) q^{\beta}$
is given by
\begin{align}\label{fourierformula}
(\theta^d f)(q) = \sum_{\beta}\left(\prod_{\sigma\in\Sigma}\sigma(\det \beta)^{d(\sigma)}\right)a(\beta) q^{\beta}.
\end{align}
\paragraph{Relationship between certain values of $\ci$ and $\padic$ Eisenstein series}
Let $R$ be an $\OK$-algebra in which $p$ splits completely, together with embeddings
\begin{align*}
i_\infty&: R\hookrightarrow\IC\\
i_p&: R\hookrightarrow R_0 := \varprojlim_n R/p^nR.
\end{align*}

As explained in \cite{kaCM} (for Hilbert modular forms) and in \cite{EDiffOps} (for automorphic forms on $U(n,n)$), the above differential operators can be used to relate certain values of $\ci$-automorphic forms to values of $p$-adic automorphic forms.  We now apply results from these papers to relate special values of the $p$-adic and $\ci$ Eisenstein series discussed earlier in this paper.  

Denote by $L$ the lattice $\OK^{2n}\subseteq K^{2n}$.  The $R$-valued points of $_{\mathcal{K}}Sh(R)$ parametrize tuples $\uA$ consisting of an abelian variety together with a polarization, endomorphism, and level structure (as discussed in \cite[Section 2.2]{EDiffOps}).  For each $R$-valued point of $_{\mathcal{K}}Sh(R)$ corresponding to an abelian variety $\uA$, denote by $\omega_{\uA/R}^{\pm}$ an ordered basis over $R$ for the sheaf of differentials $\uo_{\uA/R}^{\pm}$ (where $\uo_{\uA/R} = \uo_{\uA/R}^{+}\oplus \uo_{\uA/R}^{-}$ is the usual sheaf of one-forms on $\uA$ over $R$, as in \cite[(2.2)]{EDiffOps}).  

For each point $g\in G(\IQ)\backslash X\times G(\adeles_f)/\mathcal{K}$ (identified with the $\IC$-valued points of $_{\mathcal{K}}Sh(\IC)$ as above), let $z = g_{\infty}\gimel$, and denote by $\uA(g)$ the tuple such that $\uA(g)$ is the tuple $P_z$ (dependent upon $L$) in \cite[Section 2.3.2]{EDiffOps} consisting of a complex abelian variety with PEL structure.  Let $\uo^{\pm}(g)$ denote an ordered basis over $\IC$ for $\uo_{\uA/\IC}^{\pm}$ such that we have an equality of lattices $p_z(L) = L(\uA(g), \uo_{\uA/\IC}^{\pm}(g))\subseteq\IC^{2n}$.  The precise definitions of $P_z$, $p_z(L)$, and $L(\uA(g), \uo_{\uA/\IC}^{\pm}(g))$ are tedious and unnecessary here; they are discussed in detail in \cite[Section 2]{EDiffOps}.

\subsubsection{Comparison of certain values of $\ci$- and $\padic$-Eisenstein series}\label{cipadiccomparison}
Now suppose $\uA$ is a tuple over $R$ together with an ordered basis $\uo_{\uA/R}^{\pm}$ over $R$, such that there is an $R$-submodule $\Split(\uA/R)$ in $H^1_{DR}(\uA/R)$ inducing a splitting over $R$
\begin{align*}
\uo_{\uA/R}\oplus \Split(\uA/R)\rightarrow H^1_{DR}(\uA/R)
\end{align*}
giving an isomorphism such that $H^1_{DR}(\uA/R)^{\pm}\subseteq\uo_{\uA/R}\oplus\Split(\uo/\uA)$, simultaneously satisfying Conditions $\left(\dag\right)$ and $\left(\ddag\right)$ of \cite{EDiffOps}, namely the image of the inclusion $\Split(\uA/R)\otimes_R\IC\hookrightarrow H^1(\uA_{\IC}^{an}, \IC)$ is the the antiholomorphic differentials $H^{0, 1}$ and the image of the inclusion $\Split(\uA/R)\otimes_R R_0\hookrightarrow H^1(\uA_{R_0}^{\padic}, R_0)$ is the the unit root subspace of $H^1(\uA_{R_0}^{\padic}, R_0)$.  (i.e. $\uA$ is a CM abelian variety.)

Now, let $g\in G(\IQ)\backslash X\times G(\adeles_f)/\mathcal{K}$, and suppose that $\uA(g)$ is the extension by scalars to $\IC$ of a CM abelian variety $\uA$ over $R$ together with $\uo^{\pm}_{\uA/R}$ and $\Split(\uA/R)$ satisfying Conditions $\left(\dag\right)$ and $\left(\ddag\right)$ (explained in the previous paragraph).  Let $\Omega^{\pm}\in\gln(\IC)^{\Sigma}$ satisfy
\begin{align*}
\uo^{\pm}_{\uA/R} = \Omega^{\pm}\cdot\uo^{\pm}(g)
\end{align*}
over $\IC$, and let $c^{\pm}\in\gln(R_0)^{\Sigma}$ satisfy
\begin{align*}
\uo^{\pm}_{\uA/R} = c^{\pm}\cdot\uo_{can}^{\pm}(\uA)
\end{align*}
over $R_0$, where $\uo_{can}^{\pm}$ is the basis with this notation in \cite[Section 5]{EDiffOps}.  (For the results in the current paper, the precise definition of $\uo_{can}^{\pm}$ is unnecessary.)  Then by application of the main results on algebraicity in \cite{EDiffOps},
\begin{tiny}
\begin{align}\label{padiccx}
\det(\Omega^+)^{-(k+\nu+d)}&\det(\Omega^-)^{(\nu-d)}\det(z-{ }^t\bar{z})^{-d}\prod_{\sigma\in\Sigma}(2\pi i)^{nd(\sigma)}\Psi_{n, d, k}D(1, K, \mathfrak{b}, p, k)^{-1} E_{k+2d, \nu-d}\left(z; h, \chi, \mu, \frac{k}{2}\right)\nonumber \\
&= \det(c^+)^{-(k+\nu+d)}\det(c^-)^{(\nu-d)}\theta^d G_{k, \nu, \chi, \mu}(\uA).
\end{align}
\end{tiny}
(In \cite{EEHLS}, we explain how to choose CM periods $\Omega^{\pm}$ and $c^{\pm}$ uniformly for all CM abelian varieties at once.  For the current paper, however, this will not be necessary.)

Since we started with functions on groups, we conclude this subsection by reminding the reader of the simple relationship
\begin{align*}
E_{k+2d, \nu-d}\left(z; h, \chi, \mu, \frac{k}{2}\right) = j_{g_{\infty}}^{k+2d, \nu-d}(\gimel)E_f(m(h)g_{\infty}),
\end{align*}
where $h\in \gln(\adeles^\infty)$ and $g_{\infty}\in G(\IR)^{\Sigma}$ such that $g_{\infty}\gimel = z$.

%

\paragraph{Highest weights and higher dimensional representations}
The differential operators discussed above
are a special case of certain differential operators that map the Eisenstein series to  automorphic forms whose weight is a representation of dimension $\geq1$.  (See, for example, \cite{shar} for a description of these more general operators in the $\ci$-case.  The author's thesis, \cite{EDiffOps}, is the reference for these operators in the $p$-adic setting.) 

This special case covered in the present paper is sufficient for $p$-adically interpolating certain special values of non-holomorphic Eisenstein series $E_{k, \nu}(z; h, \chi, \mu, s)$ with $s$ not necessary $\frac{k}{2}$ and $k$ not necessarily an integer (but rather a $\Sigma$-tuple of integers).  This is enough for both \cite{EEHLS} and for the current state of homotopy theory.

In \cite{apptoSHLvv}, we use these operators to generalize the results of this paper to automorphic forms of non-scalar weights.  We obtain a similar result to Equation \eqref{weightup} but with $\Psi_{n, d, k}$ replaced by a number in terms of the highest weight of a given representation.  The results of this paper naturally generalize to the setting of \cite{apptoSHLvv}, but for simplicity, we have divided them into two portions.

\paragraph{Lie theory}
There is a Lie-theoretic approach to the $\ci$-differential operators that might be more natural in the context of automorphic forms on adele groups.  (See, for example, the appendix to \cite{shar}.)  We have used the above approach due to the fact that there are currently no references for a $p$-adic analogue of the $\ci$-Lie-theoretic differential operators.  The author hopes to fill in this gap in the literature in the near future.
%

%
%

\section{A $p$-adic Eisenstein measure}\label{EisensteinMeasureSection}

Let $T(r)$ be as in Equation \eqref{Tofr}.
Let 
\begin{align*}
\mathcal{G} = \left(\left(\OK\otimes\ZZ_p\right)^{\times}/\overline{\OK^{\times}}\right)\times T(r),
\end{align*}
 where $\overline{\OK^\times}$ denotes the closure of $\OK^\times$ in $\left(\OK\otimes\ZZ_p\right)^{\times}$.  Let $\mathcal{V}$ denote the space of $p$-adic automorphic forms, as above.  From any Hecke character $\chi,$ we obtain a $p$-adic character $\tilde{\chi}$ in the notation of Equation \eqref{tildechi} by shifting to $p$.

Recall that in Sections \ref{qexpansions} through \ref{viewedpadically}, we introduced Eisenstein series and $E_{F, \mu}$ whose $q$-expansions were shown in Equations \eqref{Gknuchimu} and \eqref{EFmu} to be
\begin{align*}
G_{k, \nu, \chi, \mu}(q) & = \sum_{\beta}\left(\alpha(\beta, \chi)\phi_{\left(\chi_1^{-1}\chi_2\mu\right)}({ }^t\beta)\prod_{v\in\Sigma}\sigma_v(\det\beta)^{k-n}\right)q^{\beta}\\
E_{F, \mu}(q) & = \sum_{\beta}\left(\left(\prod_{v\ndivides p\infty}P_{\beta, v, \mathfrak{b}}\left(F(\pi_v\right))\right)\cdot F(\det\beta^{-1})N_{E/\IQ}(\det\beta)^{-n}\phi_\mu(\beta)\right)q^{\beta}.
\end{align*}
Throughout this section, we shall consider the case in which $a(\sigma)b(\sigma) = 0$ for all $\sigma\in\Sigma$.

\begin{thm}[An Eisenstein measure]\label{measthm1}
There is a continuous $\mathcal{V}$-valued measure $d\phi$ on $\mathcal{G}$ such that
\begin{align*}
\int_\mathcal{G} (\tilde{\chi}_p, \mu\prod_{\sigma\in\Sigma}\sigma^{d(\sigma)}) d\phi = \theta^d G_{k, \nu, \chi', \mu} = E_{\tilde{\chi}, \mu\cdot\prod_{\sigma\in\Sigma}\sigma^{d(\sigma)}}
\end{align*}
for all locally constant characters $\mu$ on $T(r)$ and Hecke characters $\chi$ with infinity type $\psi_{k, \nu} = \prod_{\sigma\in\Sigma}\sigma^{-k}\left(\frac{\bar{\sigma}}{{\sigma}}\right)^{\nu}$ and conductor dividing $p^{\infty}$, where $\chi' := \chi\left|\cdot\right|_{K/\IQ}^{\frac{k}{2}}$.
\end{thm}
\begin{proof}
From the expression of the Fourier coefficients in Example \ref{Fexl} and the discussion at the end of Section \ref{certaineseries}, it immediately follows that there is a $\mathcal{V}$-valued measure $d\phi$ on $\mathcal{G}$ given by
\begin{align*}
\int_\mathcal{G} (F, \mu) d\phi = E_{F, \mu}
\end{align*}
for each continuous $\OCp$-valued function $F$ on $\left(\left(\OK\otimes\ZZ_p\right)^{\times}/\overline{\OK^{\times}}\right)$ extended by $0$ to $\OK\otimes\ZZ_p$ and each continuous $\OCp$-valued function $\mu$ on $T(r)$.  The theorem now follows directly from Equation \eqref{fourierformula} and Example \ref{Fexl} (together with the $q$-expansion principle), by taking $F = \tilde{\chi}$.
\end{proof}

As a direct corollary of Theorem \ref{measthm1}, we obtain:

\begin{cor}\label{Dpmeas}
Let $R$ be an $\OK$-subalgebra of $\OCp$.  For each point $\uA$ consisting of an ordinary abelian variety with PEL structure (i.e. a point of the Igusa tower) defined over $R$, together with a choice of basis for $\uo = \uo_{\uA/R}$, there is a $\OCp$-valued measure $d\phi(\uA, \uo)$ on $\mathcal{G}$ defined by
\begin{align*}
\int_\mathcal{G} (F, \mu) d\phi(\uA, \uo) = E_{F, \mu}(\uA, \uo)
\end{align*}
for all $\OCp$-valued continuous functions $(F, \mu)$ on $\mathcal{G}$.
\end{cor}

Technically, the above theorem and corollary were valued in the $p$-adic automorphic forms restricted to points of the Shimura variety for which $m(\prod_{v\ndivides p\infty}h_v) = 1_{2n}\in M(\adeles^{p, \infty})$.  However, the above results naturally generalize to the ring of $p$-adic automorphic forms even after we remove the restriction $m(\prod_{v\ndivides p\infty}h_v) = 1_{2n}\in M(\adeles^{p, \infty})$.
Thus, we generalize Theorem \ref{measthm1} and Corollary \ref{Dpmeas} to obtain measures on
\begin{align*}
\mathcal{G'}:=\mathcal{U}\times T(r),
\end{align*}
where $\mathcal{U}$ is the id\`ele class group $H(p^\infty) \isomto $ the maximal abelian extension of $K$ of conductor $p^{\infty}$.  We identify $\left(\left(\OK\otimes\ZZ_p\right)^{\times}/\overline{\OK^{\times}}\right)$ with its image in $\mathcal{U}$ and note that $\mathcal{U}/\left(\left(\OK\otimes\ZZ_p\right)^{\times}/\overline{\OK^{\times}}\right)$ is finite.

\begin{cor}[Improvement of Theorem \ref{measthm1}]
There is a continuous $\mathcal{V}$-valued measure $d\phi$ on $\mathcal{G'}$ such that the $q$-expansion of 
\begin{align*}
\int_\mathcal{G'} (\tilde{\chi}, \mu\prod_{\sigma\in\Sigma}\sigma^{d(\sigma)}) d\phi
\end{align*}
is $\theta^d G_{k, \nu, \chi', \mu}(\bullet, q) = E_{\tilde{\chi}, \mu\cdot\prod_{\sigma\in\Sigma}\sigma^{d(\sigma)}}(\bullet, q)$
for all locally constant characters $\mu$ on $T(r)$ and Hecke characters $\chi$ with infinity type $\psi_{k, \nu} = \prod_{\sigma\in\Sigma}\sigma^{-k}\left(\frac{\bar{\sigma}}{{\sigma}}\right)^{\nu}$ and conductor dividing $p^{\infty}$, where $\chi' := \chi\left|\cdot\right|_{K/\IQ}^{\frac{k}{2}}$.  The first variable $\bullet$ denotes a point of $M(\adeles^{p, \infty})$.
\end{cor}
\begin{proof}
Apply Lemma \ref{mavsm1} and Theorem \ref{measthm1}.
\end{proof}

\begin{thm}[$p$-adic interpolation of special values of $\ci$-Eisenstein series]
Let $d = (d(\sigma))_{\sigma\in\Sigma}\in\ZZ^{\Sigma}$, and let $R$ be an $\OK$ subalgebra of $\IC$ such that there is also an injection $R\hookrightarrow R_0 =\varprojlim_nR/p^n R$.  For each CM point $\uA(g)$ defined over $R$ corresponding to a point $g\in Sh(\IC)$ and satisfying conditions $\left(\dag\right)$ and $\left(\ddag\right)$ as in Section \ref{cipadiccomparison}, 
\begin{small}
\begin{align*}
&\frac{\int_{\mathcal{G'}} (\tilde{\chi}, \mu\cdot \prod_{\sigma\in\Sigma}\sigma^{d(\sigma)}) d\phi(\uA(g))}{\det (c^+)^{k+\nu+d}\det (c^-)^{-(\nu-d)}}  = \\
 &\frac{(2\pi i)^{nd}\Psi_{n,d,k} \det(z-{ }^t\bar{z})^{-d}}{\det(\Omega^+)^{k+\nu+d}\det(\Omega^-)^{-(\nu-d)}D(n, K, \mathfrak{b}, p, k)}E_{k+2d, \nu-d}\left(z; \prod_{v\ndivides p\infty}g_v, \chi', \mu, \frac{k}{2}\right)
 \end{align*}
\end{small}
for all locally constant characters $\mu$ on $T(r)$ and Hecke characters $\chi$ with infinity type $\psi_{k, \nu} = \prod_{\sigma\in\Sigma}\sigma^{-k}\left(\frac{\bar{\sigma}}{{\sigma}}\right)^{\nu}$ and conductor dividing $p^{\infty}$, where $\chi' := \chi\left|\cdot\right|_{K/\IQ}^{\frac{k}{2}}$ and $z = g_{\infty}\gimel$.
\end{thm}
\begin{proof}
This follows from Theorem \ref{measthm1} and its corollaries, together with Equation \eqref{padiccx}.
\end{proof}

\bibliography{eischen}   

\def\Dbar{\leavevmode\lower.6ex\hbox to 0pt{\hskip-.23ex \accent"16\hss}D}
  \def\cfac#1{\ifmmode\setbox7\hbox{$\accent"5E#1$}\else
  \setbox7\hbox{\accent"5E#1}\penalty 10000\relax\fi\raise 1\ht7
  \hbox{\lower1.15ex\hbox to 1\wd7{\hss\accent"13\hss}}\penalty 10000
  \hskip-1\wd7\penalty 10000\box7}
  \def\cftil#1{\ifmmode\setbox7\hbox{$\accent"5E#1$}\else
  \setbox7\hbox{\accent"5E#1}\penalty 10000\relax\fi\raise 1\ht7
  \hbox{\lower1.15ex\hbox to 1\wd7{\hss\accent"7E\hss}}\penalty 10000
  \hskip-1\wd7\penalty 10000\box7} \def\Dbar{\leavevmode\lower.6ex\hbox to
  0pt{\hskip-.23ex \accent"16\hss}D}
  \def\cfac#1{\ifmmode\setbox7\hbox{$\accent"5E#1$}\else
  \setbox7\hbox{\accent"5E#1}\penalty 10000\relax\fi\raise 1\ht7
  \hbox{\lower1.15ex\hbox to 1\wd7{\hss\accent"13\hss}}\penalty 10000
  \hskip-1\wd7\penalty 10000\box7}
  \def\cftil#1{\ifmmode\setbox7\hbox{$\accent"5E#1$}\else
  \setbox7\hbox{\accent"5E#1}\penalty 10000\relax\fi\raise 1\ht7
  \hbox{\lower1.15ex\hbox to 1\wd7{\hss\accent"7E\hss}}\penalty 10000
  \hskip-1\wd7\penalty 10000\box7} \def\Dbar{\leavevmode\lower.6ex\hbox to
  0pt{\hskip-.23ex \accent"16\hss}D}
  \def\cfac#1{\ifmmode\setbox7\hbox{$\accent"5E#1$}\else
  \setbox7\hbox{\accent"5E#1}\penalty 10000\relax\fi\raise 1\ht7
  \hbox{\lower1.15ex\hbox to 1\wd7{\hss\accent"13\hss}}\penalty 10000
  \hskip-1\wd7\penalty 10000\box7}
  \def\cftil#1{\ifmmode\setbox7\hbox{$\accent"5E#1$}\else
  \setbox7\hbox{\accent"5E#1}\penalty 10000\relax\fi\raise 1\ht7
  \hbox{\lower1.15ex\hbox to 1\wd7{\hss\accent"7E\hss}}\penalty 10000
  \hskip-1\wd7\penalty 10000\box7} \def\cprime{$'$}
\providecommand{\bysame}{\leavevmode\hbox to3em{\hrulefill}\thinspace}
\providecommand{\MR}{\relax\ifhmode\unskip\space\fi MR }
\providecommand{\MRhref}[2]{%
  \href{http://www.ams.org/mathscinet-getitem?mr=#1}{#2}
}
\providecommand{\href}[2]{#2}
\begin{thebibliography}{AHR10}

\bibitem[AHR10]{AHR}
Matthew Ando, Michael Hopkins, and Charles Rezk, \emph{Multiplicative
  orientations of {$KO$}-theory and of the spectrum of topological modular
  forms}, 2010, \url{http://www.math.uiuc.edu/~mando/papers/koandtmf.pdf}.

\bibitem[Beh09]{beh}
Mark Behrens, \emph{Eisenstein orientation}, 2009, Typed notes available at
  \url{http://www-math.mit.edu/~mbehrens/other/coredump.pdf}.

\bibitem[EHLS]{EEHLS}
Ellen~E. Eischen, Michael Harris, Jian-Shu Li, and Christopher~M. Skinner,
  \emph{$p$-adic {$L$}-functions for unitary groups}, In preparation.

\bibitem[Eis]{apptoSHLvv}
Ellen~E. Eischen, \emph{A higher-dimensional $p$-adic {E}isenstein measure for
  unitary groups}, In preparation.

\bibitem[Eis12]{EDiffOps}
\bysame, \emph{$p$-adic differential operators on automorphic forms on unitary
  groups}, Annales de l'Institut Fourier \textbf{62} (2012), no.~1, 177--243.

\bibitem[Har81]{hasv}
Michael Harris, \emph{Special values of zeta functions attached to {S}iegel
  modular forms}, Ann. Sci. \'Ecole Norm. Sup. (4) \textbf{14} (1981), no.~1,
  77--120. \MR{MR618732 (82m:10046)}

\bibitem[Har86]{ha86}
\bysame, \emph{Arithmetic vector bundles and automorphic forms on {S}himura
  varieties. {II}}, Compositio Math. \textbf{60} (1986), no.~3, 323--378.
  \MR{MR869106 (88e:11047)}

\bibitem[HLS06]{SHL}
Michael Harris, Jian-Shu Li, and Christopher~M. Skinner, \emph{{$p$}-adic
  {$L$}-functions for unitary {S}himura varieties. {I}. {C}onstruction of the
  {E}isenstein measure}, Doc. Math. (2006), no.~Extra Vol., 393--464
  (electronic). \MR{MR2290594 (2008d:11042)}

\bibitem[Hop95]{hopkins94}
Michael~J. Hopkins, \emph{Topological modular forms, the {W}itten genus, and
  the theorem of the cube}, Proceedings of the {I}nternational {C}ongress of
  {M}athematicians, {V}ol.\ 1, 2 ({Z}\"urich, 1994) (Basel), Birkh\"auser,
  1995, pp.~554--565. \MR{1403956 (97i:11043)}

\bibitem[Hop02]{hopkinsICM}
M.~J. Hopkins, \emph{Algebraic topology and modular forms}, Proceedings of the
  {I}nternational {C}ongress of {M}athematicians, {V}ol. {I} ({B}eijing, 2002)
  (Beijing), Higher Ed. Press, 2002, pp.~291--317. \MR{1989190 (2004g:11032)}

\bibitem[Kat78]{kaCM}
Nicholas~M. Katz, \emph{{$p$}-adic {$L$}-functions for {CM} fields}, Invent.
  Math. \textbf{49} (1978), no.~3, 199--297. \MR{MR513095 (80h:10039)}

\bibitem[Lan08]{la}
Kai-Wen Lan, \emph{Arithmetic compactifications of {PEL}-type {S}himura
  varieties}, 2008, Ph.D. thesis, Harvard University, available at
  \url{http://www.math.princeton.edu/~klan/articles/cpt-PEL-type-thesis-single%
.pdf}.

\bibitem[Lan11]{lanalgan}
\bysame, \emph{Comparison between analytic and algebraic constructions of
  toroidal compactifications of {PEL}-type {S}himura varieties}, Crelle's
  Journal (2011), To appear.

\bibitem[Shi82]{shconfluent}
Goro Shimura, \emph{Confluent hypergeometric functions on tube domains}, Math.
  Ann. \textbf{260} (1982), no.~3, 269--302. \MR{669297 (84f:32040)}

\bibitem[Shi83]{sheseries}
\bysame, \emph{On {E}isenstein series}, Duke Math. J. \textbf{50} (1983),
  no.~2, 417--476. \MR{705034 (84k:10019)}

\bibitem[Shi84]{shclassical}
\bysame, \emph{On differential operators attached to certain representations of
  classical groups}, Invent. Math. \textbf{77} (1984), no.~3, 463--488.
  \MR{759261 (86c:11034)}

\bibitem[Shi97]{sh}
\bysame, \emph{Euler products and {E}isenstein series}, CBMS Regional
  Conference Series in Mathematics, vol.~93, Published for the Conference Board
  of the Mathematical Sciences, Washington, DC, 1997. \MR{MR1450866
  (98h:11057)}

\bibitem[Shi00]{shar}
\bysame, \emph{Arithmeticity in the theory of automorphic forms}, Mathematical
  Surveys and Monographs, vol.~82, American Mathematical Society, Providence,
  RI, 2000. \MR{MR1780262 (2001k:11086)}

\bibitem[Tan99]{tan}
Victor Tan, \emph{Poles of {S}iegel {E}isenstein series on {${\rm U}(n,n)$}},
  Canad. J. Math. \textbf{51} (1999), no.~1, 164--175. \MR{1692899
  (2000e:11073)}

\end{thebibliography}

\end{document}